# On a computational approach for the approximate dynamics of averaged variables in nonlinear ODE systems: toward the derivation of constitutive laws of the rate type

Amit Acharya* and Aarti Sawant

Civil and Environmental Engineering, Carnegie Mellon University, Pittsburgh, PA 15213, U.S.A

**Abstract**

A non-perturbative approach to the time-averaging of nonlinear, autonomous ODE systems is developed based on invariant manifold methodology. The method is implemented computationally and applied to model problems arising in the mechanics of solids.

## 1. Introduction

We are interested in exploring the possibility of developing a completely-defined (i.e. closed) rate-type evolution equation for time averages of prescribed functions of state of an underlying physical system, the latter assumed to be adequately described by a system of Ordinary Differential Equations (ODE). We refer to the time averages of the prescribed state functions as coarse variables and the degrees of freedom of the underlying ODE system as fine variables. Of particular interest are the following questions: can such a closed description be generated at all? If so, is such a description generally adequate for probing the evolution of the coarse quantities when appended by an initial condition on the coarse variables alone?

It is shown in this paper that the answer to the first question is in the affirmative, both conceptually and computationally, with significant practical benefits. As for the second question, the answer is found to be negative in a strict sense. However, from a practical point of view, the proposed method can be shown to be no worse than the philosophy of physical experimentation where one might commit oneself to a choice of coarse variables and the existence of an autonomous, rate-type constitutive equation for these variables and proceed to determine the form of the right-hand-side of such an equation by experimentally monitoring small increments in response out of prescribed coarse states, *assuming* that a coarse state corresponds to a unique increment out of that state. In fact, a major component of the long-term, practical objective of this work is to work with experimentalists in making a choice of coarse variables such that over a broad range of experimental trials the incremental response of the chosen set of coarse variables is found to be unique, once a coarse state has been specified. Of course, our method lays bare the fact that such uniqueness ought not be expected in general, if one were to choose the number and type of coarse variables simply on physical grounds and convenience of measurement. Indeed, a simplification of results of Inertial Manifold Theory (Temam, 1990) to finite dimensional systems suggests that an autonomous coarse response (corresponding to the class of dissipative fine dynamics that is within the rigorous purview of Inertial Manifold Theory) is guaranteed only if the number of coarse variables adopted is bounded below by the dimension of the attractor of the fine dynamics. The fact that approximate, consistent, lower-dimensional (than the dimension of the attractor) coarse dynamics can be developed in principle and demonstrated in practice in the context of small systems displaying complicated dynamics is the subject of our prior work (Acharya 2005; Sawant and Acharya 2005), where connections of our approach to related work in the literature are presented. Such consistent coarse dynamics, however, requires information on initial conditions more than what is prescribed merely by a knowledge of coarse initial states and displays response characteristic of non-autonomous systems in general, features related to the prediction of memory effects solely as a geometric consequence of a drastic reduction in degrees of freedom in transitioning from fine to coarse dynamics. In this paper, 'memory' in dynamics refers to the possibility of different rates of evolution of a set of evolving variables from an identical state for this set.

* Corresponding author: Tel. (412) 268 4566; Fax. (412) 268 7813; email: acharyaamit@cmu.edu

Our overall approach is based on the use of invariant manifolds of the fine dynamics as proposed in Muncaster (1983), with conceptually and practically non-trivial extensions made in Acharya (2005) and Sawant and Acharya (2005). However, these approaches as well as coarse-graining approaches based on Inertial Manifold Theory (Foias et al., 1988) Center Manifold Theory (Carr and Muncaster, 1983 a, b; Roberts, 2003), or the Constructive Methods of Gorban and Karlin (Gorban et al., 2004) have not explored the possibility of running *time-averages* of fine response as coarse variables, perhaps for the reason that such a choice does not fit in naturally into the invariant manifold methodology for model reduction. This is so because such averages cannot be defined as functions of the instantaneous fine state. In this paper, we propose a simple augmentation of the fine dynamics that allows such a choice of coarse variables, followed by the use of invariant manifold methodology for obtaining closed dynamics for the time averages. To our knowledge, such an approach for dealing with time-averaging is new.

Obtaining the effective dynamics of nonlinear ODE systems is a vast field with many approaches, as summarized in the excellent review of Givon et al. (2004). Primary approaches other than invariant manifold methodology have been the method of averaging for Hamiltonian systems (e.g. Neishtadt, 2004; the weak convergence method, Bornemann, 1998) and the projection operator technique (e.g. Chorin et al., 2000). Our primary interest being in the mechanics of solids and materials science, we have been strongly influenced by the results of Abeyaratne et al. (1996), especially in their demonstration of the importance of the averaging of microscopic dynamics, as opposed to the use of averaged energetics, in the construction of effective dynamics. Our main goal is the development of computational methodology based on sound, even if not mathematically rigorous, ideas that can be put to use in the development of constitutive response in the mechanics of solids. We refer to our overall approach as the method of Parametrized Locally Invariant Manifolds (PLIM); the first word in the name is meant to emphasize the practical importance of the parametrization of locally invariant manifolds of the fine dynamics and, in particular, the evolution of the parameters along fine trajectories on these manifolds. The second word is inserted to emphasize the crucial importance of allowing for, at least conceptually, the presence of many locally invariant manifolds in fine phase space corresponding to each local domain of coarse phase space, as developed in Acharya (2005) and Sawant and Acharya (2005).

The conceptual and methodological basis of our multiscale/coarse-graining method is compared with others, both theoretical and computational, in Acharya (2005) and Sawant and Acharya (2005). The primary feature of our computational idea that differentiates it from other approaches for multiscale numerics is the suggestion that an appropriate collection of parametrizations of locally invariant manifolds actually be computed and stored ‘off-line’, a calculation that does not require evolution of the fine system. With this database in hand, a closed coarse dynamics is defined that is *consistent* with the fine (microscopic) dynamics in a well-defined sense, a property that is not transparent and/or valid for all multiscale computational approaches. Depending on the choice of coarse variables, a user-input in the procedure, the coarse theory is computable much more efficiently than the fine theory. As we show in our past work and here, this strategy allows us to compute the coarse response of complicated nonlinear systems without special treatment and regardless of whether there exists a ‘separation of scale’ – e.g. see Sawant and Acharya (2005) for an application to the Lorenz system with a choice of coarse variables that neither induces nor exploits any separation of scales. Of course, our methodology is equally capable of inducing a separation of scale by an appropriate choice of coarse variables and exploiting such a separation as we show in the examples of this paper, thus providing one avenue for computationally dealing with nonlinear, stiff systems, even oscillatory. As is well known, the ease afforded by the existence of separation of scales in a problem for doing theory does not necessarily transfer to the design of computational algorithms for such, especially for nonlinear stiff problems.

Importantly, our overall approach (Acharya, 2005; Sawant and Acharya, 2005; this paper) implies that one may not expect an autonomous coarse dynamics simply because the coarse variables are chosen as ‘slow’ and/or ‘aggregate’ variables and a separation of scales exists in this limited, literal, sense. Even more importantly, it suggests a means to achieve, in principle, a consistent coarse dynamics in the absence of a separation of scales in the sense that the evolution of the fast and slow variables do not decouple, ideas that can be turned into practice approximately. The latter feature is in contrast to some coarse-graining exercises for ODE, e.g. the homogenization of the associated linear, transport PDE as shown in Menon (2002), where the non-uniqueness of the evolution of weak limits (coarse variables) from a specified coarse state is recovered, but an understanding of the failure at the level of suggesting possible deterministic remedies is not. While recourse to stochastic models may be the simplest alternative to get some acceptable answers in the absence of a true separation of scales, it must be realized that, at a fundamental level, the non-autonomous content – both memory of coarse variables and 'noise' - of the appropriate stochastic model for coarse-graining a deterministic, nonlinear system of ODE is intimately linked to understanding of the 'orthogonal dynamics' of the eliminated variables and its coupling to the dynamics of the retained coarse variables, as laid bare by the Mori-Zwanzig Projection Operator Technique (Chorin et al., 2002).

This paper is organized as follows: the basic theoretical ideas behind our method are presented in Section 2. Section 3 contains a demonstration of the general strategy through model problems. Scalar and vector gradient flow dynamics with 'wiggly' energies due to Abeyaratne et al. (1996) and Menon (2002) are homogenized in time by our method. While low-dimensional, these are complicated problems that do not satisfy any ergodicity hypothesis for the fast dynamics, when suitably phrased as problems of averaging. Next, the macroscopic 'stress-strain' response of an atomic chain based on the Frenkel-Kontorova model (e.g., Nabarro, 1987) is computed under monotonic loading. In solving these problems, it is assumed that full knowledge of fine initial conditions is available. However, in most applications it is reasonable to assume only some limited knowledge, at best, of the possible fine states corresponding to a coarse initial condition. In Section 4, *a* strategy is proposed to deal with such situations within our general approach, based purely on considerations of practicality. The paper ends with some concluding remarks in Section 5.

## 2. A simple technique for time-averaging

The autonomous fine dynamics is defined as

$$\begin{aligned} \frac{df}{dt}(t) &= H\left(f(t)\right) \\ f(0) &= f_*. \end{aligned} \tag{1}$$

$f$ is an $N$-dimensional vector of fine degrees of freedom and $H$ is a generally nonlinear function of fine states, denoted as the vector field of the fine dynamical system. Equation $(1)_2$ represents the specification of initial conditions. $N$ can be large in principle, and the function $H$ rapidly oscillating.

Let $\Lambda$ be a user-specified function of the fine states producing vectors with $m$ components whose time averages over intervals of period $\tau$ can be measured in principle and are of physical interest. Given the fixed time interval $\tau$ characterizing the resolution of coarse measurements in time, a coarse trajectory corresponding to each fine trajectory $f(\cdot)$ is defined as the following running time average:

$$c(t) = \frac{1}{\tau}\int_t^{t+\tau} \Lambda\left(f(s)\right) ds. \tag{2}$$

Roughly speaking, it is a closed statement of evolution for $c$ that we seek. The statement is unambiguous only after we specify what sort of initial conditions we may want to specify. For the purpose of this section we assume that fine initial conditions are known with certainty. Then, the goal is to develop a closed evolution equation for $c$, i.e. an equation that can be used for evolving $c$ without concurrently evolving (1), corresponding to fine trajectories out of a prescribed set of fine initial conditions.

Clearly,

$$\frac{dc}{dt}(t) = \frac{1}{\tau}\left[\Lambda\left(f(t+\tau)\right) - \Lambda\left(f(t)\right)\right]. \tag{3}$$

If we now introduce a forward trajectory $f_f(\cdot)$ corresponding to a trajectory $f(\cdot)$ as

$$f_f(t) := f(t+\tau), \tag{4}$$

then

$$\frac{df_f}{dt}(t) = H\left(f_f(t)\right). \tag{5}$$

Also, given an initial state $f_*$ we denote by $f_{**}$ the state defined as the solution of (1) evaluated at time $\tau$. With these definitions in hand, we augment the fine dynamics (1) to

$$\begin{aligned}
&\frac{df_f}{dt}(t) = H\left(f_f(t)\right)\\
&\frac{df}{dt}(t) = H\left(f(t)\right)\\
&f_f(0) = f_{**}\\
&f(0) = f_*
\end{aligned} \tag{6}$$

and apply invariant manifold techniques to (6). In detail, on an $m$- dimensional coarse phase space whose generic element we denote as $c$, we seek functions $G_f$ and $G$ that satisfy the first-order, quasilinear partial differential equations

$$\left.\begin{aligned}\sum_{k=1}^{m}\frac{\partial G_f^I}{\partial c^k}\left(\frac{1}{\tau}\left[\Lambda^k\left(G_f\right)-\Lambda^k\left(G\right)\right]\right)=H^I\left(G_f\right)\\ \sum_{k=1}^{m}\frac{\partial G^I}{\partial c^k}\left(\frac{1}{\tau}\left[\Lambda^k\left(G_f\right)-\Lambda^k\left(G\right)\right]\right)=H^I\left(G\right)\end{aligned}\right\}\quad I=1\text{ to }N \tag{7}$$

at least locally in $c$-space. Assuming that we have such a pair of functions over the domain containing the point

$$c_* := c(0) \tag{8}$$

defined from (1) and (2) which, moreover, satisfies the conditions

$$\begin{aligned}G_f\left(c_*\right)&=f_{**}\\ G\left(c_*\right)&=f_*,\end{aligned} \tag{9}$$

it is easy to see that a local-in-time fine trajectory defined by

$$\begin{aligned}\Gamma_f\left(t\right)&:=G_f\left(c\left(t\right)\right)\\ \Gamma\left(t\right)&:=G\left(c\left(t\right)\right)\end{aligned} \tag{10}$$

through the coarse local trajectory satisfying

$$\begin{aligned}&\frac{dc}{dt}=\frac{1}{\tau}\left[\Lambda\left(G_f\left(c\right)\right)-\Lambda\left(G\left(c\right)\right)\right]\\ &c\left(0\right)=c_*\end{aligned} \tag{11}$$

is the solution of (6) (locally). A solution pair $\left(G_f, G\right)$ of (7) represents a parametrization of a locally invariant manifold of the dynamics (6). By a locally invariant manifold we mean a set of points in phase space such that the vector field of (6) is tangent to the set at all points. Thus, a trajectory of (6) enters and exits a locally invariant manifold only through the boundary of the manifold.

Also, note that if $\left(\widehat{G}_f, \widehat{G}\right)$ and $\left(\breve{G}_f, \breve{G}\right)$ are two solutions to (7) and (9) on an identical local domain in $c$-space containing $c_*$ and $\widehat{c}\left(\cdot\right)$ and $\breve{c}\left(\cdot\right)$ are the corresponding coarse trajectories defined as solutions to (11), then local uniqueness of solutions to (6) implies

$$\begin{aligned}\widehat{G}_f\left(\widehat{c}\left(t\right)\right)&=:\widehat{\Gamma}_f\left(t\right)=\breve{\Gamma}_f\left(t\right):=\breve{G}_f\left(\breve{c}\left(t\right)\right)\\ \widehat{G}\left(\widehat{c}\left(t\right)\right)&=:\widehat{\Gamma}\left(t\right)=\breve{\Gamma}\left(t\right):=\breve{G}\left(\breve{c}\left(t\right)\right).\end{aligned} \tag{12}$$

Thus,

$$\begin{aligned}&\frac{d\widehat{c}}{dt}\left(t\right)=\frac{d\breve{c}}{dt}\left(t\right)\quad\text{locally in time,}\\ &\widehat{c}\left(0\right)=\breve{c}\left(0\right)\end{aligned} \tag{13}$$

from (11), and assuming $\widehat{c}$ and $\breve{c}$ are continuous, $\widehat{c}\equiv\breve{c}$, locally.

Hence, given *any* pair of mappings $G_f, G$ satisfying (7) *and* (9) on a domain containing $c_*$, we consider (11) as the *consistent, closed* theory for the evolution of the coarse variables $c$. Obstruction to the construction of solutions to (7) is explored in Acharya (2005), providing one reason for seeking multiple local solutions as implemented in Sawant and Acharya (2005).

We mention here that the choice of coarse variables (2) can be appended with functions of the instantaneous fine state of (1) to define more coarse variables. The determining equations for the $G$ functions and the coarse evolution for these variables are different in form from (7) and (11). This class of coarse variables is the more standard choice in applications of invariant manifold methodology as outlined, e.g., in Muncaster (1983). Applied loads are introduced as such coarse variables in the model problems of this paper in Section 3.

The basic algorithmic idea (Sawant and Acharya, 2005), then, is as follows: consider a region of the phase-space of (1), trajectories/trajectory-segments lying in which one is interested in coarse-graining. By definition of the forward-shifted variables, this delineates a region of the phase-space of (6), say $R$, within which trajectory segments of (6) are to be coarse-grained. One constructs as many pairs of functions $\left(G_f, G\right)$ locally in $c$-space as possible, with continuous extensions to the boundaries of their local domains of definition in coarse space. The ideal limit is when one can ensure that each point of $R$ belongs to the range of at least one such pair. Given a fine initial condition for (6), a choice of a locally invariant manifold is made such that this initial condition belongs to it. The coarse theory (11) is now evolved. Generally, the coarse trajectory evolves to a coarse state, say $c_b$, on the boundary of the local coarse domain of definition of the picked pair $\left(G_f, G\right)$, currently in use in evolving (11). In the ideal case, another pair of solutions to (7) should now be available such that $\left(G_f\left(c_b\right), G\left(c_b\right)\right)$ belongs to its range, and with the choice of this pair, the coarse evolution is further continued. In practice, only a finite number of local solutions of (7) can be computed and not all points of $R$ can be included in the range of some of these solutions. In such a case, on reaching the boundary of a local coarse domain using the pair $\left(G_f, G\right)$, one evaluates the fine state

$$\left(G_f\left(c_b\right) + \Delta t\, H\left(G_f\left(c_b\right)\right), G\left(c_b\right) + \Delta t\, H\left(G\left(c_b\right)\right)\right), \tag{14}$$

where $\Delta t$ is the time-step of the fine evolution. Suppose that this fine state does not belong to the range of any of the precomputed finite collection of local solutions to (7). Then, one merely chooses the nearest (in some reasonable metric) available manifold to this fine state and carries on with the coarse calculation.

In this sense, every choice of a precomputed set of locally invariant manifolds defines an approximate coarse dynamical system. As more and more manifolds are computed, approaching the ideal limit defined in the previous paragraph, the approximate coarse dynamical system becomes exact. While practically somewhat useless in circumstances involving a large number of fine degrees of freedom, this conceptual exercise is an important consistency guarantee for our coarse-graining scheme attesting to its correctness.

The conceptual importance of the condition (9) for consistent coarse-graining in the presence of a mathematically arbitrary choice of coarse variables, motivated solely by physical considerations, cannot be overestimated. As the considerations leading up to (13) show, it would not be possible to guarantee a correct coarse-graining by simply choosing any solution of (7) that does not satisfy (9). This raises the unsettling, but inevitable, prospect that the evolution of a small set of physically realistic coarse variables of autonomous fine dynamics may not be posed correctly as initial value problems in the traditional sense whose trajectories are uniquely defined by simply specifying coarse initial conditions, and that fine initial conditions might exert some influence on coarse dynamics. On the other hand, practical experience with modeling of inelastic response in solid mechanics shows that rate-type coarse equations for a relatively small number of variables representing the coarsening of an enormous collection of fine entities are often practically adequate for describing macroscopic behavior. These two observations taken together lead us to propose a strategy in Section 4 for making a choice for the computation of a distinguished set of parametrizations for locally invariant manifolds that enable the development of an approximate, rate-type coarse dynamics which, however, would display memory effects for reasons demonstrated in Sawant and Acharya (2005).

## 3. Model Problems

In this section we illustrate through examples that the method can be implemented successfully. The first two examples belong to the category where the coarse variables are running time averages of the fine variables. We choose two model problems whose dynamics is of the gradient flow type governed by energy functions with small perturbations. The first example describes the kinetics of a phase transforming material (Abeyaratne et al., 1996) and the second example is a gradient flow in 2 dimensions (Menon, 2002).

The third example belongs to the category where the coarse variables are representative of a space-time average of a collection of fine degrees of freedom. Specifically, we model the behavior of a 1-d, discrete chain of atoms in a nonlinear force field, with linear nearest-neighbor interactions.

In all the examples presented in this section the fine system is evolved numerically, and the time-averaged response obtained from the fine system is denoted as the *actual*/*actual coarse* response. The response obtained by evolving the coarse theory is denoted as *coarse* response.

Both the fine and the coarse theories are evolved discretely in time with identical numerical schemes for all the reported comparisons.

### 3.1 Kinetics of a phase transforming material with wiggly energy

This example is selected from the paper by Abeyaratne et al. (1996), where the evolution of volume fraction of one of the variants of martensite in a plate subjected to biaxial loading was studied. A fine scale gradient flow model is developed and subsequently averaged to obtain a nonstandard macroscopic kinetic law in closed form through a weak convergence argument. Due to the restriction of our method to autonomous fine systems, we modify the cyclic loading program of Abeyaratne et al. (1996) slightly to obtain the following fine system of equations describing a gradient flow type kinetic law for the volume fraction $\lambda$ , under a loading program characterized by stresses $(\sigma_1, \sigma_2, \sigma_3)$:

$$\begin{aligned} &\dot{\lambda} = -\mu \frac{dW}{d\lambda}(\lambda, \sigma_1, \sigma_2), \lambda(0) = \lambda_0, \\ &\dot{\sigma}_1 = \eta\,\sigma_3, \dot{\sigma}_3 = \eta(\sigma_0 - \sigma_1), \sigma_2(t) = \text{constant} \\ &\sigma_1(0) = \sigma_1^0; \sigma_3(0) = \sigma_3^0. \end{aligned} \tag{15}$$

Here $\sigma_1, \sigma_2$ represent uniform applied tractions in 1 and 2 directions along the specimen boundary. The variable $\eta$ represents the frequency of the periodic loading function and $\sigma_3$ is an auxiliary variable introduced to allow an oscillatory response in $\sigma_1$. Figure 1(a) shows the periodically varying stress $\sigma_1$ and constant stress $\sigma_2$, along with one of the loading programs in Abeyaratne et al. (1996) that is approximated by $(15)_{3,4,5}$.

The term $\mu$ represents the mobility of the system and $\lambda_0$ is the initial volume fraction. The term $W$ is the total energy of the specimen comprising the terms

$$\begin{aligned} W_{load} &= -\begin{Bmatrix} \lambda^2(\sigma_1^2 + \sigma_2^2)(\alpha^2 - \gamma^2)^2 \big/ (\alpha^2 + \gamma^2) \\ +2\lambda(\sigma_1^2\gamma^2 - \sigma_2^2\alpha^2)(\alpha^2 - \gamma^2) \big/ (\alpha^2 + \gamma^2) \\ +(\sigma_1\gamma + \sigma_2\alpha)^2 \end{Bmatrix}^{\frac{1}{2}} \\ W_{trans.layer} &= c_1\lambda^2 + c_2(1 - \lambda^2) \\ W_{pert} &= a\varepsilon\cos\left(\frac{\lambda}{\varepsilon}\right) \ ; \ \varepsilon \to 0 \\ W_0 &= W_{load} + W_{trans.layer} \\ W &= W_0 + W_{pert}, \end{aligned} \tag{16}$$

where $\alpha, \gamma, c_1, c_2, a, \varepsilon$ are parameters. For the physical development of the total energy function, consult Abeyaratne et al. (1996).

The main observation of their paper related to the effect of the third type of energy $W_{pert}$ on the dynamics. This term physically characterizes the tip-splitting phenomenon of twinned bands within laminates under the loading. Here $a$ is the intensity of the perturbation and $\varepsilon \ll 1$ is a small parameter characterizing the wavelength of the energy wiggles. Figure 1(b) shows the change in energy due to the contribution of the energy $W_{pert}$. The crucial observation due to Abeyaratne et al. (1996) is that this small addition to $W_0$, that averages out to zero, completely changes the force governing the dynamics of the problem, as seen in Figure 1(c, d), i.e. a small contribution to energy does not necessarily constitute a small contribution to the driving force as can be directly seen by taking a derivative of $(16)_3$ with respect to $\lambda$. Under the periodic loading program, the volume fraction $\lambda$ shows a fixed or frozen response under varying load in certain regions of loading and this is attributed (Abeyaratne et al., 1996) to the situation where $\lambda$ gets stuck in a local minimum of $W$, as can also be inferred from a plot of the forces 1(d).

We study the evolution of volume fraction by using $W$ as the energy function with $\varepsilon \ll 1$ as a fixed parameter. It is important to note that the development of the coarse theory in our approach is non-asymptotic in the sense that it does not depend on taking the limit $\varepsilon \to 0$, a feature that has both merits and demerits in the context of homogenization in time.

**Case 1: Cyclic loading Program:**

Following the procedure in Section 2, we define a time averaged volume fraction $\overline{\lambda}$ using the definition (2), with $\Lambda : f \mapsto f$. The corresponding augmented variable as defined in (4) is $\lambda_f$. The fine theory (15) is augmented with an evolution equation corresponding to $\lambda_f$ of the form

$$\dot{\lambda}_f(t) = -\mu \frac{dW}{d\lambda}\left(\lambda_f(t), \sigma_1(t+\tau), \sigma_2(t+\tau)\right), \ \lambda_f(0) = \lambda_{f0}. \tag{17}$$

Figure 2(a) shows the evolution of volume fraction $\lambda$, averaged $\overline{\lambda}$, and the time shifted variable $\lambda_f$ in time. Figure 2(b) shows the comparison of the actual and averaged volume fractions $\lambda, \overline{\lambda}$ as a function of the difference in applied stresses $\sigma_1 - \sigma_2$. For different choices of $\tau$ (time interval used for averaging), the averaged volume fraction $\left(\overline{\lambda}\right)$ profiles (obviously) differ from the fine $(\lambda)$ trajectory, but contain essential features of fine response (frozen $\overline{\lambda}$ or flat portions seen in Figure 2).

We choose the averaged volume fraction $\overline{\lambda}$ and the stresses $\sigma_1, \sigma_3$ of the loading program as the coarse variables. Because the solution for $\sigma_1, \sigma_3$ in (15) can be expressed explicitly in terms of trigonometric functions, it is possible to write down the time –shifted stress terms appearing in (17) in terms of $\sigma_1(t), \sigma_3(t)$ and functions of the parameter $\tau$. Thus the fine system is not augmented for these variables. We denote the driving force for $\lambda_f$ obtained by this rewriting as $\widehat{(\partial W/\partial \lambda)}$, which is a function of $\lambda_f, \sigma_1, \sigma_3$ and the parameter $\tau$. As the stress $\sigma_2$ is constant, we treat it as a fixed parameter while computing the locally invariant manifolds and evolving the coarse dynamics.

We denote the fine variables by functions $\lambda = G_1(\bar{\lambda}, \sigma_1, \sigma_3)$ and. The solutions to $G_1, G_2$ are obtained by solving the governing equations (18) obtained from the augmented fine theory:

$$\begin{aligned} &\frac{\partial G_2}{\partial \bar{\lambda}}\left(\frac{G_2 - G_1}{\tau}\right) + \frac{\partial G_2}{\partial \sigma_1}\eta\sigma_3 + \frac{\partial G_2}{\partial \sigma_3}\eta(\sigma_0 - \sigma_1) + \mu\widehat{\left(\frac{\partial W}{\partial \lambda}\right)}(G_2, \sigma_1, \sigma_2; \tau) = 0 \\ &\frac{\partial G_1}{\partial \bar{\lambda}}\left(\frac{G_2 - G_1}{\tau}\right) + \frac{\partial G_1}{\partial \sigma_1}\eta\sigma_3 + \frac{\partial G_1}{\partial \sigma_3}\eta(\sigma_0 - \sigma_1) + \mu\frac{\partial W}{\partial \lambda}(G_1, \sigma_1, \sigma_2) = 0. \end{aligned} \tag{18}$$

Here the first term in $(18)_{1,2}$ corresponds to the left hand side of $(7)_{1,2}$, with $m = 1$. The next two terms arise due to the dependence of the functions $G_1, G_2$ on the other coarse variables $\sigma_1, \sigma_3$. The equation (18) is solved by using the Least squares finite element method and EI-TC method described in Sawant (2005) and the Appendix of this paper[i].

With the solutions obtained for (18), the coarse theory corresponding to equation (11) can be written as,

$$\begin{aligned} &\dot{\bar{\lambda}} = \frac{1}{\tau}\left(G_2(\bar{\lambda}, \sigma_1, \sigma_3) - G_1(\bar{\lambda}, \sigma_1, \sigma_3)\right) \\ &\dot{\sigma}_1 = \eta\sigma_3, \quad \dot{\sigma}_3 = \eta(\sigma_0 - \sigma_1). \end{aligned} \tag{19}$$

To evolve the coarse theory we obtain a consistent initial condition $\bar{\lambda}(0)$ by evolving the fine system up to time $\tau$. The initial conditions $\sigma_1(0), \sigma_3(0)$ are selected as in (15).

To illustrate the advantage of using the PLIM method, we compute the coarse response from numerical integration of the equations (19) by using much larger time steps than the maximal time step required to obtain the response of the fine system (15). The maximal time-step for the fine evolution is defined as that value of the time-step such that any time-step smaller than the maximal one reproduces the same computed response as the maximal one. The latter are subsequently $\tau$ -time averaged according to (2) for comparison with coarse response. We refer to the $\tau$ -averaged response of a computed fine solution as the *actual coarse* response. The ratio of the coarse to fine steps is denoted as $c/f$. We present the coarse response computed with $c/f$ ratios of 1, 10, and 100.

Figures 3, 4 show the comparisons of coarse response obtained from the coarse and fine theory with parameter values $\alpha = 1.0619$, $\gamma = 1.0231$, $c_1 = 0.017$ MPa, $c_2 = 0.0255$ MPa, $a = 0.025$ MPa, $\mu$=5.4 MPa$^{-1}$s$^{-1}$ taken from Abeyaratne et al. (1996). The small parameter $\varepsilon$ and the time scale $\tau$ are selected as $0.005$ and $2.0\,\text{sec}\,(200 \times \text{maximal fine time-step})$ respectively.

---

[i] We thank Dr. Hao Huang for discussions on this matter.

The initial condition on $\sigma_1, \sigma_3$ are selected as $\sigma_1(0) = \sigma_o + \sigma_d$ and $\sigma_3(0) = 0$, with a mean stress $\sigma_0 = 5.81\,MPa$ and maximum amplitude of $\sigma_3$ as $\sigma_d = \pm 1\,MPa$; the stress $\sigma_2 = 5.85\,MPa$ is held constant. The time period of the sinusoidal loading program is selected as $20\sec$. The initial volume fraction was assigned as $\lambda(0) = 0.8$. For this problem, the averaging time-scale parameter $\tau$ turns out to be approximately 5 times the period of oscillation observed in the evolution of the fine volume fraction and $1/10^{th}$ of the period of loading program. The consistent initial conditions for $\lambda_f(0)$ and $\bar{\lambda}(0)$ are obtained by evolving the fine system (15) and (17) up to time $\tau$. Figure 3 shows the comparison of actual and coarse response obtained by using the coarse-to-fine time-step ratio as $c/f = 1$. Figure 3(a) shows the comparison of fine $\lambda(\cdot)$ and $G_1(\bar{\lambda}(\cdot), \sigma_1(\cdot), \sigma_3(\cdot))$ solutions retrieved from the locally invariant manifolds corresponding to the coarse trajectory $(\bar{\lambda}(\cdot), \sigma_1(\cdot), \sigma_3(\cdot))$ obtained by evolving the coarse theory. Figure 3(b) shows the comparison between the actual coarse solution, as defined above, and the coarse theory solution for averaged volume fraction $\bar{\lambda}$ as a function of the difference in the applied stresses $\sigma_1 - \sigma_2$.

Figure 3 shows that the coarse response is in good agreement with the actual response, and the computed manifolds capture the regions with the wiggles and the flat regions in the volume fraction accurately. Figure 4 shows the coarse response computed with two different ratios of coarse-to-fine $(c/f)$ time-steps. Figures 4(a), 4(b) correspond to coarse response with $c/f = 10$ and Figures 4(c), 4(d) correspond to coarse response with $c/f = 100$.

In the examples above we use a finite set of locally invariant manifolds that predicts reasonable coarse response through the coarse theory. In general, the volume fraction $\lambda$ settles into a limit cycle after some initial transient, where the actual response is periodic. Our coarse theory solutions match well with the actual coarse response for fine initial conditions on the limit cycle, but for initial conditions off of the limit cycle some initial discrepancies are seen even though the limit cycle behavior is recovered. Figure 5 shows an example of this behavior. Figure 5(a) shows the actual and coarse evolution of $\bar{\lambda}$ in time. Figure 5(b) shows the comparison of actual and coarse averaged volume fraction $\bar{\lambda}$ as a function of the difference in the applied stresses $\sigma_1 - \sigma_2$ (hysteresis). The selected initial condition in this case is off of the limit cycle. After some initial transient the coarse trajectory is attracted towards the limit cycle and the coarse response shows good agreement with the actual response.

**Case 2: Creep Test:**

For this case, the specimen is subjected to uniform loading tractions $\sigma_1, \sigma_2$, which are held constant during the evolution of the volume fraction. The fine system for this case is also (15) and (17) with $\eta = 0, \sigma_3^0 = 0$. We choose $\bar{\lambda}$ and $\sigma_1$ as the coarse variables with $\sigma_2$ treated as a fixed parameter. Of course, $\sigma_1$ could also be treated as a parameter in this problem, just as $\sigma_2$ could have been treated as a coarse variable in both problems with $\dot{\sigma}_2 \equiv 0$. The fine variables are denoted as $\lambda = G_1\left(\bar{\lambda}, \sigma_1\right)$ and $\lambda_f = G_2\left(\bar{\lambda}, \sigma_1\right)$, which represent 2D locally invariant manifold. The governing equation for $G_1, G_2$ are

$$\begin{aligned} &\frac{\partial G_2}{\partial \bar{\lambda}}\left(\frac{G_2 - G_1}{\tau}\right) + \mu \frac{\partial W}{\partial \lambda}\left(G_2, \sigma_1, \sigma_2\right) = 0 \\ &\frac{\partial G_1}{\partial \bar{\lambda}}\left(\frac{G_2 - G_1}{\tau}\right) + \mu \frac{\partial W}{\partial \lambda}\left(G_1, \sigma_1, \sigma_2\right) = 0. \end{aligned} \tag{20}$$

We mention in passing that we calculate manifolds for the two loading programs separately only as a device to simplify matters in these 'proof of principle' demonstrations. In general, assuming all the loading programs can be subsumed within one general autonomous ODE description, the variables and parameters of the loading description can be included as independent variables for the calculation of manifolds. The parameters of the loading description can be treated as both fine and coarse variables with no evolution, this being a standard trick of invariant manifold methodology.

With the solutions to (20) in hand, the coarse theory takes the form

$$\dot{\bar{\lambda}} = \frac{1}{\tau}\left(G_2\left(\bar{\lambda}, \sigma_1\right) - G_1\left(\bar{\lambda}, \sigma_1\right)\right). \tag{21}$$

Figure 6 shows three examples with different initial conditions $\lambda(0)$ and the applied stress $\sigma_1$. For this case the applied stresses $\sigma_1, \sigma_2$ are held constant and the coarse theory (21) is evolved till the volume fraction is completely frozen. The initial condition for $\lambda_f, \bar{\lambda}$ are obtained by evolving the fine system (15) and (17) up to time $\tau$ ($2$ sec). The Figures 6(a), 6(b) show the evolution of volume fractions in time, from the three different initial conditions for $c/f = 1$. Figure 6(a) shows the actual $\lambda(.)$ and $G_1\left(\bar{\lambda}(.), \sigma_1(.)\right)$ similar to Figures 3, 4. Figure 6(b) shows the comparison of actual and coarse averaged volume fraction $\bar{\lambda}$ as a function of time. The coarse response obtained for this case also matches well with the exact response, and also with the results presented in Abeyaratne et al (1996).

Figure 7 shows the comparison of volume fractions $\lambda, \bar{\lambda}$ as a function of time, with the coarse response computed with $c/f$ ratios of 10 and 100. Figures 7(a), 7(b) correspond to coarse response with $c/f = 10$, which are in excellent agreement with the actual response. The Figures 7(c), 7(d) correspond to coarse response with $c/f = 100$, and here the coarse theory solution shows numerical error due to the larger time steps used. However it shows the correct trend in the solution.

### *3.2 2D gradient system with wiggly energy*

This example is selected from the paper by Menon (2002), where a 2-d gradient system is studied with a governing energy function of the following form

$$
\begin{aligned}
&W_o = \frac{1}{2}\left(\lambda_1 y^2 + \lambda_2 z^2\right), \\
&W_{pert} = \left(R + r\cos\left(\frac{z}{\varepsilon}\right)\right)\sin\left(\frac{y}{\varepsilon}\right).\cos(\beta) + r\sin\left(\frac{z}{\varepsilon}\right)\sin(\beta) \quad \varepsilon \to 0.
\end{aligned} \tag{22}
$$

We use the parameter values

$$R = 2, r = 1, \beta = \pi/3, \lambda_1 = 1, \lambda_2 = 1 \tag{23}$$

as one representative set suggested in Menon (2002). The small parameter $\varepsilon \ll 1$ is kept fixed. The fine system for this example takes the form

$$
\begin{aligned}
&\dot{y} = -\frac{\partial\left(W_o + W_{pert}\right)}{\partial y} = H_y(y,z), \quad \dot{z} = -\frac{\partial\left(W_o + W_{pert}\right)}{\partial z} = H_z(y,z) \\
&y(0) = y_0, z(0) = z_0.
\end{aligned} \tag{24}
$$

Numerical calculation of trajectories of the fine system indicate that this system has 3 distinct parts to its fine phase space shown in Figure 8, qualitatively consistent with the analytical conclusions of Menon (2002). Figure 8 shows two lines $\rho = \infty, \rho = 0$ (dashed), corresponding to the ratio $\rho = \dot{z}/\dot{y}$. These lines divide the fine phase space in zone I and zone II. The shaded region is denoted as zone III. The Figure shows different trajectories (solid) starting in the zone I, which contains no fixed points. The trajectories, approach zone III (containing fixed points) where the trajectories show sticking only. In this region the movement of trajectory is almost negligible. In zone II (transition zone), the dynamics is a combination of sticking and slipping in that trajectories evolve a little bit when initiated in this region before coming to a halt.

The coarse variables are selected as running time averages of the variables $y, z$ (i.e. $(c_1, c_2) = (\overline{y}, \overline{z})$ *and* $\Lambda : f \mapsto f$ in (2)). Following the procedure in section 2, two new fine variables $y_f, z_f$ are introduced that are similar to $\lambda_f$ in the previous example. The fine system (24) is augmented with the evolution equations for the variables $y_f, z_f$ of the form

$$\dot{y}_f = H_y\left(y_f, z_f\right), \quad \dot{z}_f = H_z\left(y_f, z_f\right), \quad y_f(0) = y_{f0}, \quad z_f(0) = z_{f0}. \tag{25}$$

The variables of the augmented fine system are denoted as $y = G_1(\overline{y}, \overline{z}), y_f = G_2(\overline{y}, \overline{z})$ and $z = G_3(\overline{y}, \overline{z}),\ z_f = G_4(\overline{y}, \overline{z})$. Thus we obtain a set of 4 functions representing locally invariant manifolds by solving the equations

$$\begin{aligned}
&\frac{\partial G_1}{\partial \overline{y}}\left(\frac{G_2-G_1}{\tau}\right)+\frac{\partial G_1}{\partial \overline{z}}\left(\frac{G_4-G_3}{\tau}\right)=H_y\left(G_1,G_3\right)\\
&\frac{\partial G_2}{\partial \overline{y}}\left(\frac{G_2-G_1}{\tau}\right)+\frac{\partial G_2}{\partial \overline{z}}\left(\frac{G_4-G_3}{\tau}\right)=H_y\left(G_2,G_4\right)\\
&\frac{\partial G_3}{\partial \overline{y}}\left(\frac{G_2-G_1}{\tau}\right)+\frac{\partial G_3}{\partial \overline{z}}\left(\frac{G_4-G_3}{\tau}\right)=H_z\left(G_1,G_3\right)\\
&\frac{\partial G_4}{\partial \overline{y}}\left(\frac{G_2-G_1}{\tau}\right)+\frac{\partial G_4}{\partial \overline{z}}\left(\frac{G_4-G_3}{\tau}\right)=H_z\left(G_2,G_4\right).
\end{aligned} \tag{26}$$

In terms of these computed set of $G$ functions, the coarse theory takes the form

$$\begin{aligned}
\dot{\overline{y}}&=\frac{1}{\tau}\left(G_2\left(\overline{y},\overline{z}\right)-G_1\left(\overline{y},\overline{z}\right)\right)\\
\dot{\overline{z}}&=\frac{1}{\tau}\left(G_4\left(\overline{y},\overline{z}\right)-G_3\left(\overline{y},\overline{z}\right)\right).
\end{aligned} \tag{27}$$

To evolve the coarse theory (27) we obtain consistent initial conditions $\overline{y}\left(0\right),\overline{z}\left(0\right)$ by evolving the fine system (24), (25) up to time $\tau$ from fine initial conditions.

For this model problem, two cases with different values of $\varepsilon$ are presented. In the following figures a comparison of actual and coarse response for 5 trajectories with different coarse-to-fine time-step ratios is presented. Figure 9 shows the trajectories in coarse phase space $\overline{y}-\overline{z}$ with the choice $\varepsilon=0.01$ and $\tau=0.2\,\text{sec}$, whereas Figure 10 shows the same with the choice $\varepsilon=0.001$ and $\tau=0.02\,\text{sec}$. (Again, for both cases $\tau$ is $200$ times the maximal fine time-step). The Figure also shows the fixed zone $\left(zone-III\right)$. In both the Figures it is seen that for the various $c/f$ ratios, the coarse response is in agreement with the actual coarse response in $zone-I$, which is a region with no fixed points. The coarse response diverges from the exact one in the transition $zone-II$, which shows stick-slip behavior and the movement of trajectories is almost negligible once they reach the boundary of the fixed zone (zone $-III$).

Interestingly, in this problem it is observed that most fine trajectories in the fine phase space can be enveloped into only one surface described by the coarse variables $\overline{y},\overline{z}$ for the choice of parameters listed in (23).

### *3.3 Macroscopic stress-strain curve of an atomic chain*

Figure 11 shows a chain of atoms placed on a substrate that exerts a spatially periodic force on the atoms derived from a potential of period $a_s$. The neighboring atoms are connected by linear elastic spring (or segments) with stiffness $k$, and $a_o$ is the distance between the atoms when the springs have zero strain. The chain is fixed at one end and an external load is applied at the free end. Before application of the load the all the atoms are placed in the troughs of the substrate potential such that $a_s = a_o$. Thus, the chain is in a stationary state corresponding to an absolute minimum of potential energy. We refer to this state as a *zero-strain state* in rest of this section.

The corresponding mechanical model introduced by Frenkel and Kontorova (e.g. Braun, 1949), can be derived from the standard Hamiltonian,

$$h = \frac{m_a}{2}\sum_i \left(\frac{dx_i}{dt}\right)^2 + \frac{p_s}{2}\sum_i \left(1 - \cos\left(\frac{2\pi\left(x_i - x_1\right)}{a_s}\right)\right) + \frac{k}{2}\sum_i \left(x_{i+1} - x_i - a_0\right)^2 \quad i = 2\ to\ \eta + 1 \tag{28}$$

where $\eta + 1$ is the total number of particles(atoms). The first term is the kinetic energy, where $m_a$ is the mass of each particle and $x_i$ is the position of the $i^{th}$ particle in the chain. The second term is a part of the potential energy, which characterizes the interaction of the chain with an external periodic substrate potential, where $p_s, a_s$ are the amplitude and period of the potential respectively. The last potential energy term takes into account a linear coupling between the nearest neighbors of the chain, and is characterized by the elastic constant $k$ and the equilibrium distance of the inter-particle potential, $a_o$, in the absence of the substrate potential.

From (28), the equations of motion of the particles in the chain subjected to external load are obtained as

$$m_a \frac{d^2 x_i}{dt^2} + \frac{\pi p_s}{a_s}\sin\left(\frac{2\pi\left(x_i - x_1\right)}{a_s}\right) - k\left(x_{i+1} - 2x_i + x_{i-1}\right) = F_i, \tag{29}$$

where $F_i$ represents the external load applied at the $i^{th}$ particle. We consider a chain with $a_s = a_o$ (zero-strain state) as shown in Figure 11, which helps to replace the position $x_i$ with the displacements $u_i$ using $x_i = (i-1)a_s + u_i$. The normalized version of (29), in terms of displacements, is

$$\frac{d^2 \hat{u}_i}{d\hat{t}^2} + \sin\left(2\pi\hat{u}_i\right) - \hat{k}\left(\hat{u}_{i+1} - 2\hat{u}_i + \hat{u}_{i-1}\right) = \hat{F}_i. \tag{30}$$

Here, all the forces in (29) are normalized with respect to ratio of the magnitude of the substrate potential $\varepsilon_s$ and its period $a_s$ i.e. $\hat{F} = F/\left(\pi p_s/a_s\right)$. The actual displacement $u_i$ is normalized as $\hat{u} = u/a_s$; the stiffness is normalized as $\hat{k} = k/\left(\pi p_s/a_s^2\right)$ and the time is normalized as $\hat{t} = t/\left(\sqrt{m_a a_s^2/\pi p_s}\right)$. *For clarity in the rest of the equations the carets are dropped from the non-dimensional quantities*. We introduce the velocities $v_i$ and a term $\alpha_{ij}$ representing the elastic force in equation (29)-(30). Thus the normalized fine ODE system that we work with is written as follows:

$$\left.\begin{aligned} \dot{u}_i &= v_i \\ \dot{v}_i &= -\sin\left(2\pi u_i\right) + \alpha_{ij}\, u_j + F_i \end{aligned}\right\} \quad i = 2\ to\ \eta + 1$$
$$u_1 = 0, v_1 = 0\ at\ \ fixed\ end; \quad F_i = 0\ \ i \neq \eta + 1; \quad F_{\eta+1} = L(t), \dot{L} = f_r. \tag{31}$$

where $L(t)$ is the normalized external load, applied to $\eta+1^{th}$ atom at the end of the chain. We restrict our study to the case of monotonically increasing or decreasing load, i.e. $L(t) = f_o + f_r t$, where $f_r$ is a loading rate.

We consider the coarse variable to be a space and time averaged strain. The time averaged displacement and velocity $\bar{u},\bar{v}$, are defined as

$$\begin{gathered}\bar{u}_i(t) = \frac{1}{\tau}\int_t^{t+\tau} u_i(s)\,ds, \quad \bar{v}_i(t) = \frac{1}{\tau}\int_t^{t+\tau} v_i(s)\,ds, \\ \bar{\varepsilon}_i(t) = \bar{u}_{i+1}(t) - \bar{u}_i(t),\end{gathered} \tag{32}$$

where $\bar{\varepsilon}_i$ is the time-averaged strain in $i^{th}$ spring (or segment) coupling the adjacent atoms $i, i+1$. We define the space and time averaged strain as the number average of time-averaged strains in all $\eta$ couplings,

$$\tilde{\varepsilon}(t) = \frac{1}{\eta}\sum_{i=1}^{\eta}\bar{\varepsilon}_i(t) = \frac{\bar{u}_{\eta+1}(t) - \bar{u}_1(t)}{\eta} \tag{33}$$

Here, $\bar{u}_{\eta+1}$ is a time averaged displacement of the $\eta+1^{th}$ particle at the end of the chain where the load is applied.

We choose the space-time averaged strain $\tilde{\varepsilon}$ and the external force $L$ as the coarse variables (i.e. $c_1 = \tilde{\varepsilon}, \Lambda = \int_0^L u_x dx$ in (2)) and treat the loading rate $f_r$ as a fixed parameter. Using the definition of space-time averaged strain $\tilde{\varepsilon}$ in (33), we obtain an evolution equation for it as

$$\dot{\tilde{\varepsilon}}(t) = \frac{\dot{\bar{u}}_{\eta+1}}{\eta} = \frac{\left(u_{\eta+1}(t+\tau) - u_{\eta+1}(t)\right)}{\tau\eta}. \tag{34}$$

and $(31)_5$ is the evolution equation for $L$.

Following the recipe set forth in Section 2, the fine system is augmented with time shifted variables $u_f, v_f$ corresponding to each particle. The corresponding evolution equations are

$$\begin{gathered}\left(\dot{u}_f\right)_i = \left(v_f\right)_i, \quad \left(\dot{v}_f\right)_i = -\sin\left(\left(2\pi u_f\right)_i\right) + \alpha_{ij}\left(u_f\right)_j + \hat{F}_i, \quad i = 2 \; to \; \eta+1 \\ \hat{F}_i = 0 \quad i \neq \eta+1, \quad \hat{F}_{\eta+1} = L(t+\tau) = L(t) + \tau f_r, \\ \left(u_f\right)_1 = 0, \left(v_f\right)_1 = 0 \quad at \;\; fixed \; end.\end{gathered} \tag{35}$$

where $\hat{F}_i$ is the time-shift in the external load applied at $i^{th}$ particle. It should be noted that we do not introduce an augmented variable corresponding to the external load, as it can be expressed as a function of the coarse variable $L$ and the time shift $\tau$ $(35)_4$.

Now the fine displacements and velocities for each atom are denoted as $u_{k+1} = G_k(\tilde{\varepsilon}, L)$, $\left(u_f\right)_{k+1} = G_{k+\eta}(\tilde{\varepsilon}, L)$, $v_{k+1} = G_{k+2\eta}(\tilde{\varepsilon}, L)$, and $\left(v_f\right)_{k+1} = G_{k+3\eta}(\tilde{\varepsilon}, L)$, where $k = 1 \, to \, \eta$. (Note that we do not define the functions $G$ for the degrees of freedom at the fixed end).

With these definitions, a set of $4\eta$ functions are computed by solving the following governing equations:

$$\left.\begin{aligned}
&\frac{\partial G_k}{\partial \tilde{\varepsilon}}\left(\frac{G_{2\eta}-G_\eta}{\tau}\right)+\frac{\partial G_k}{\partial L} f_r = G_{k+2\eta},\\
&\frac{\partial G_{k+\eta}}{\partial \tilde{\varepsilon}}\left(\frac{G_{2\eta}-G_\eta}{\tau}\right)+\frac{\partial G_{k+\eta}}{\partial L} f_r = G_{k+3\eta},\\
&\frac{\partial G_{k+2\eta}}{\partial \tilde{\varepsilon}}\left(\frac{G_{2\eta}-G_\eta}{\tau}\right)+\frac{\partial G_{k+2\eta}}{\partial L} f_r = -\sin\left(2\pi G_k\right)+\alpha_{kl}\, G_l + F_k,\\
&\frac{\partial G_{k+3\eta}}{\partial \tilde{\varepsilon}}\left(\frac{G_{2\eta}-G_\eta}{\tau}\right)+\frac{\partial G_{k+3\eta}}{\partial L} f_r = -\sin\left(2\pi G_{k+\eta}\right)+\alpha_{kl}\, G_{l+\eta} + \hat{F}_k.
\end{aligned}\right\} k = 1\, to\, \eta. \tag{36}$$

A collection of locally invariant manifolds, represented by the functions $G$, is computed for different values of loading rate $f_r$. Using this information, the coarse theory takes the form,

$$\begin{aligned}
\dot{\tilde{\varepsilon}} &= \frac{1}{\tau\,\eta}\left(G_{2\eta}\left(\tilde{\varepsilon},L\right)-G_\eta\left(\tilde{\varepsilon},L\right)\right)\\
\dot{L} &= f_r.
\end{aligned} \tag{37}$$

In the computed examples below, we consider a chain of 11 atoms (i.e. $\eta = 10$), subjected to three loadings of the type $L(t) = f_o + f_r t$ corresponding to different values of $f_0$ and $f_r$. In all the examples, the initial conditions on displacement and velocity are zero, which correspond to the static *zero-strain state* described in Figure 11.

The coarse response is obtained by evolving the coarse theory (37) and a consistent initial condition $\tilde{\varepsilon}(0)$ for the coarse variable is obtained by evolving the fine system(31), (35) up to a time $\tau$. The value of $\tau$ is selected to be 200 times the maximal time-step required for fine evolution, the latter as defined in Section 3.1. Here $\tau$ is approximately $1/3^{rd}$ ($1/5^{th}$) of the period of the highest (lowest) mode of free vibration of the linearized system corresponding to (29) about the zero-strain state (i.e. $\sin\left(2\pi u_i\right)\approx 2\pi u_i, F_i = 0$ in (29)). It should be kept in mind that the stiffness of the linearized system changes drastically as the base state for linearization samples the states encountered during nonlinear evolution.

Figure 12 show two examples, where the atomic chain is subjected to a monotonically increasing load, with $f_o = 0$ and loading rate $f_r$ as 0.6 and 0.8 respectively[ii] **.**

---

[ii] The non-dimensional loading rates approximately corresponds to $0.25\, N/\text{sec}$ and the time interval $\tau$ is approximately of the order of $10^{-11}$ sec . These rough estimates of $f_r, \tau$ are obtained by using a lumped mass system of an axial steel bar. Here the stiffness of connecting linear springs is equal to the axial stiffness and the amplitude of the substrate potential and applied load are selected in terms of the stiffness and initial length of the bar.

Figure 12 shows the variation of the space-time averaged strain $\tilde{\varepsilon}$ as a function of monotonic loading function $L$, obtained by evolving the fine system (31), (35) along with (33). Figure 12(a) shows the coarse and actual coarse responses representing the solutions obtained from coarse and fine theory respectively. The plot also shows the coarse response corresponding to the linear system without the substrate potential. The coarse responses in this figure are computed by evolving (37) using the same time step required for the fine system $(c/f = 1)$. It is interesting to observe the similarities of the coarse response with macroscopic stress-strain curves, even though the model is extremely simplified. Figure 12(b) is a magnified version of Figure 12(a) which shows behavior similar to that in the linear elastic regime and stage II- stage III transition in work-hardening response for elastic plastic solids. For the cases with initial configuration of *zero-strain state*, the force-strain curve is symmetric in tension and compression as shown in Figure 13.

Figure 14 shows the coarse response computed by evolving the coarse theory with larger time steps, with coarse to fine time steps ratios $(c/f)$ of 10 and 100. This plot also shows good agreement between coarse and actual coarse results for all $(c/f)$ ratios. Figure 15 shows the coarse response for a third example, where the chain is loaded impulsively with $f_o = 1.5$ and subsequently the load is increased monotonically at a constant loading rate $f_r = 0.8$. This plot shows the coarse responses for $(c/f)$ ratios of 1, 10 and 100. The force-strain curve is similar to the work-hardening region, as the applied impulsive load is much larger than the yielding loads.

**4. Strategy for defining an approximate rate-type coarse dynamics: a practical expedient**

The considerations in Section 2 indicate that the ideal collection of parametrizations (mappings) of locally invariant manifolds is one where each point of the fine phase-space region of interest is an element of the range of at least one parametrization in the collection. For a fine dynamics with a large number of degrees of freedom with a choice of coarse variables that are few in number, attaining this goal is impractical. It is also true that empirical, memory-less, rate-type constitutive theories for the nonlinear, inelastic behavior of solids have been reasonably successful as modeling tools. In addition, it should be kept in mind that procedures like those of Enskog, Chapman, Grad, and Truesdell and Muncaster[iii] in the Kinetic Theory of Gases, the prototypical example of coarse-graining a nonlinear dynamical systems with many degrees of freedom, commit to the complete determination of a strict rate-type coarse theory that would conceptually parallel developing (11) with *only one* local solution of (7). Such a choice may be physically interpreted as the *assumption* that the evolution of coarse states is not sensitive to fine

[iii] The reader is referred to Truesdell and Muncaster (1980) for an excellent exposition of all of these procedures.

initial conditions, while low-dimensional dynamics indicate that this assumption cannot be valid in general (e.g. Sawant and Acharya, 2005). Finally, in many practical situations exact knowledge of fine initial conditions consistent with an observable coarse one may simply not be available.

The above reasons motivate us to propose the following compromise for the calculation of a set of parametrizations for use in the coarse theory. Arbitrarily fix a coarse state $c$, considered as a representative coarse initial condition. Denote by $\Im_c$ the set of fine states consistent with this coarse initial condition:

$$\Im_c = \left\{ f_0 : c = \int_0^\tau \Lambda\left( f(s) \right) ds;\ \frac{df}{dt} = H(f);\ f(0) = f_0 \right\}. \tag{38}$$

For each $c$, let a probability density function, $\rho(f;c)$, be defined on the set $\Im_c$ (with obvious extension to the whole of the fine phase space, $\Phi$, of (1)) that reflects the probability of a fine state being the actual fine initial condition realized, corresponding to the observed coarse initial condition $c$. Information of a physical nature may be utilized here. For example, the Maximum Entropy Postulate defining a set of quasi-equilibrium states (Gorban et al.; 2001, 2004) along with approximate theoretical information on the existence of exponentially attracting globally invariant manifolds of the fine dynamical system may be used.

In the absence of any prior information on evolution, we associate the projection $f^c$ (in some suitable norm) of the fine state

$$\int_\Phi f \rho(f;c) df \tag{39}$$

on the set $\Im_c$ as the approximate fine initial condition corresponding to the coarse initial condition $c$. In case the Maximum Entropy postulate is used to prescribe the initial probability density distribution, it is important to note that the present proposal hinges in no way on the primary constrained minimization problem involved in the implementation of the Maximum Entropy postulate having a unique solution or on the set of fine states $\left\{ f^c \right\}$, obtained as $c$ varies over the entire coarse phase space, forming a manifold.

In a local coarse region containing $c$, one now solves (7) and (9) with the identifications

$$\begin{aligned} f_* &= f^c \\ f_{**} &= f(\tau);\ \frac{df}{dt} = H(f);\ f(0) = f^c. \end{aligned} \tag{40}$$

In terms of an algorithm to implement the above idea, one picks a finite set of representative points in the coarse domain and carries out the above procedure corresponding to these points to pre-compute and store a set of locally invariant manifolds for use in the coarse theory, as indicated in Section 2.

Note that the above procedure allows a coarse state to be associated with more than one $(G_f, G)$ mapping, and the use of these in the coarse theory can result in different coarse evolution out of the single coarse state in question. If one has no prior information on how a coarse state $c$ is attained by the system, then one uses the $(G_f, G)$ described in this Section for the initial coarse evolution. On the other hand, if one attains the same state $c$ in the course of coarse evolution from the same or some other coarse initial condition through the approximate coarse dynamical system, then one has some prior information on how such a state is achieved and it is only natural that the evolution out of this state can be different in this case from when the state is treated as an initial condition for evolution.

We realize that the procedure suggested here for dealing with uncertainty in the knowledge of initial data is different from the question of transporting the initial probability density according to the fine dynamics and asking for state quantities at subsequent times as appropriate averages with respect to the evolving probability density function. While we see obvious pros and cons for both approaches, a thorough comparative evaluation awaits further research.

## 5. Concluding Remarks

Regardless of the coarse-graining technique employed, it is perhaps fair to say that given a fine dynamics containing $N$ degrees of freedom, a closed, dynamics for $M$ $(M \le N)$ coarse variables, requiring initial data on only these $M$ variables and valid for all $t \to \infty$, is possible only if there exists an $M$-dimensional globally invariant manifold in the $N$-dimensional phase space of the fine system. In essence, our method shows that an $m$-dimensional $(m < M)$ closed, coarse dynamics can also be defined, but one which requires initial conditions on at least $M$ variables to be exact. Such an exact $m$-dimensional dynamics displays memory effects in coarse variables. These general conclusions on the qualitative nature of exact, reduced coarse dynamics, i.e. memory and requirement of fine initial conditions, are similar to those obtained from the Mori-Zwanzig Projection Operator Technique (e.g. Chorin et al., 2000) of non-equilibrium statistical physics, but here reached through completely different logical arguments of a simple geometric nature. We take satisfaction in the fact that our entire procedure can be understood with only a background in multivariable calculus, the elementary theory of ordinary differential equations and equally elementary background in the theory of first-order, quasilinear, partial differential equations. The crux of the conceptual implementation of the method involves 'filling up' an $M$ dimensional manifold with $m$-dimensional, locally invariant submanifolds, in the sense of coverings. A useful geometric picture is to think of a region of three-dimensional space being filled by an infinite collection of 2-dimensional, possibly intersecting, surface patches; almost all (1-dimensional) aperiodic fine trajectories initiated in the $3$-dimensional region run along these patches, jumping appropriately from one to another, twisting and turning to cover the $3$-dimensional region densely.

The practical implementation of our method involves, obviously, the computation of a finite collection of $m$ -dimensional, locally invariant submanifolds and consequently, some error is to be expected. But the value of the fact that the source of this error is transparent in the methodology suggesting, at least, brute-force improvements to it in the form of the pre-computation and storage of a more extensive set of locally invariant manifolds, is not to be underestimated. These ideas are demonstrated in the context of the Lorenz system in Sawant and Acharya (2005). Moreover, it is not completely clear whether the validity of the coarse-grained dynamics for the limit $t \to \infty$ is necessary for computational coarse graining schemes so that the necessity of the existence of $M$ -dimensional fine, *globally* invariant manifolds may become irrelevant and one might be able to gainfully apply our technique for parametrizing particular bounded regions of fine phase space that fine trajectories to be coarse grained could very well exit. The thrust of our work in this paper is to show that our strategy can actually be executed in quite difficult nonlinear problems involving homogenization in time. Our results also show that the proposed approach can serve to set up approximate models of coarse behavior at least in the neighborhood of a collection of specific fine trajectories that may be deemed to be especially important on physical grounds, as in the example of the atomic chain.

As for error control in practical problems, we appeal to making a good choice of the coarse variables, in particular, on the knowledge that the chosen coarse variables are physical, macroscopic observables that more-or-less evolve in a state-dependent way, i.e. come close to forming a set capable of parametrizing a physically relevant invariant manifold of the fine dynamics. The main task then, is to derive an approximate right-hand-side of the coarse evolution equation. We emphasize here that our modest goal is simply the determination of the right-hand-side, *given* the definition of the coarse variables along fine trajectories. An equally valid question is understanding the process of making this choice of coarse variables with a view to making the best possible one, but as one of the most important results from Inertial Manifold Theory shows, this can lead to questions of estimating the geometry, in particular, the fractal dimension, of complicated sets (attractors) in fine phase-space that sets a lower bound on the number of coarse variables.

Running space-time averaging of PDE systems (e.g. Acharya and Roy (2005) for an example related to dislocation mechanics in crystalline materials) naturally furnish examples of the type of coarse variables we seek. With some physically motivated approximations related to the time dependence of boundary conditions on the averaging domain, such averaging problems can be fitted into our framework. Of course, there is no guarantee in such a case that the adopted coarse variables ought to display more-or-less memory-less behavior and that the approximations related to the boundary conditions are appropriate, and a future challenge for us is to apply, evaluate, and improve the methodology developed herein for such problems. In the case of a transition from dislocation mechanics to plasticity, it would also be interesting to compare the resulting dynamics to the autonomous dynamics deduced by Puglisi and Truskinovsky (2005) via approximate homogenization.

## 6. Acknowledgment

Support for this work from the Program in Computational Mechanics of the US ONR (N00014-02-1-0194) and the US AFOSR (F49620-03-1-0254) is gratefully acknowledged.

**List of Figures**

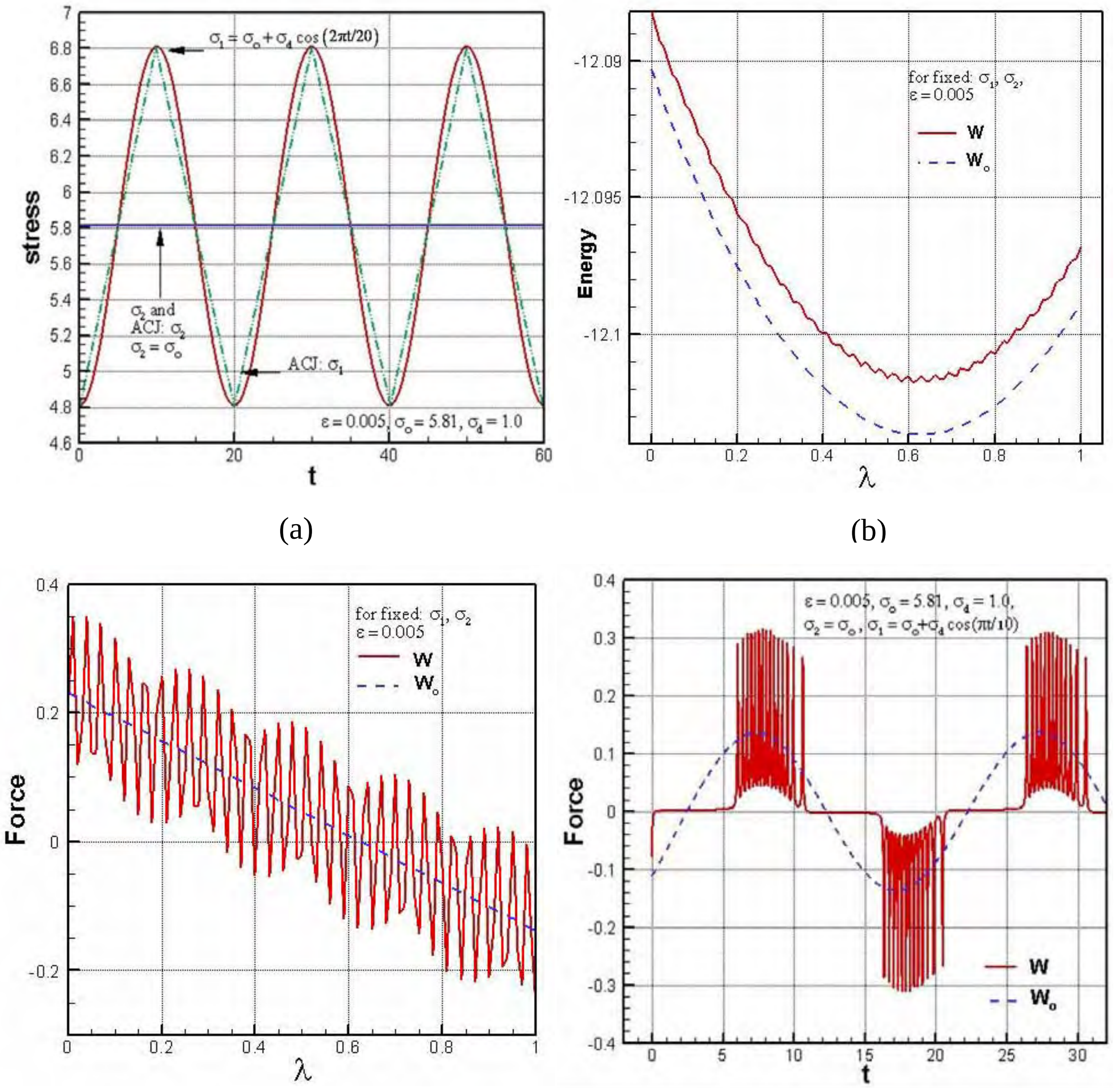

(a) (b)

(c) (d)

Figure 1: (a) Loading program $\sigma$ *vs.* $t$, (b) Energy *vs.* $\lambda$,
(c) Force *vs.* $\lambda$, (d) Force *vs.* $t$

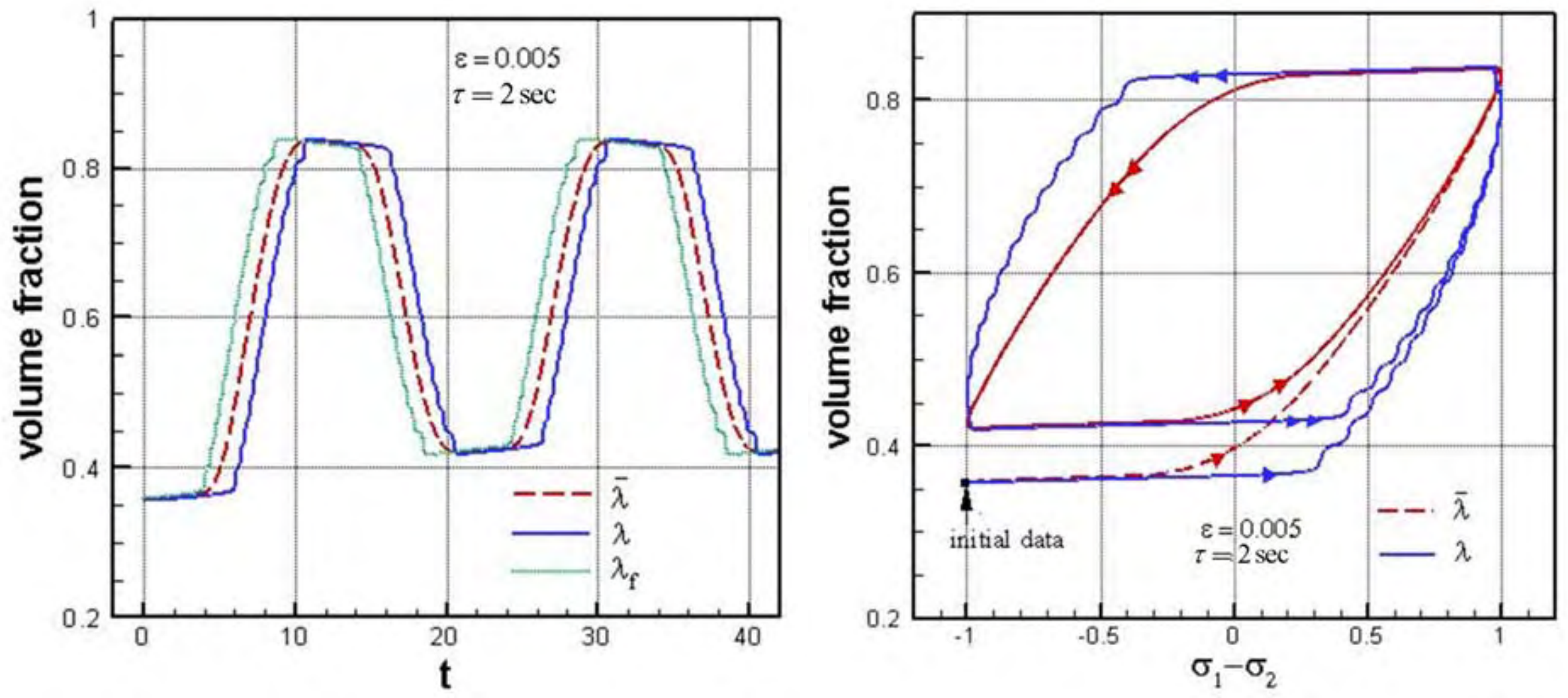

(a) (b)

Figure 2: (a) Evolution of volume fraction in time
(b) Hysteresis in volume fraction

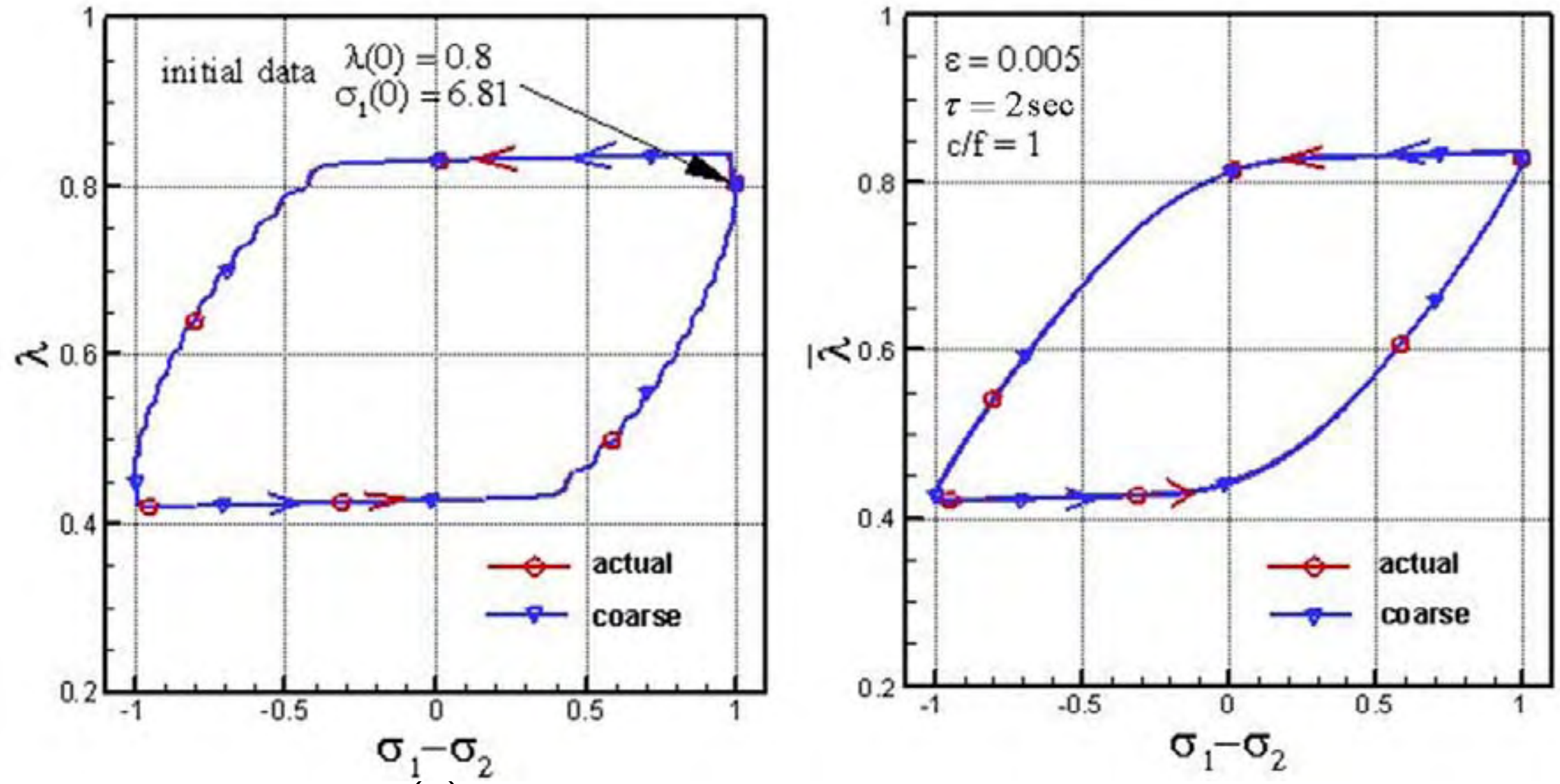

(a) (b)

Figure 3: Hysteresis. (a) fine $\lambda\, vs.\, \sigma_1-\sigma_2$ (b) averaged $\bar{\lambda}\, vs.\, \sigma_1-\sigma_2$

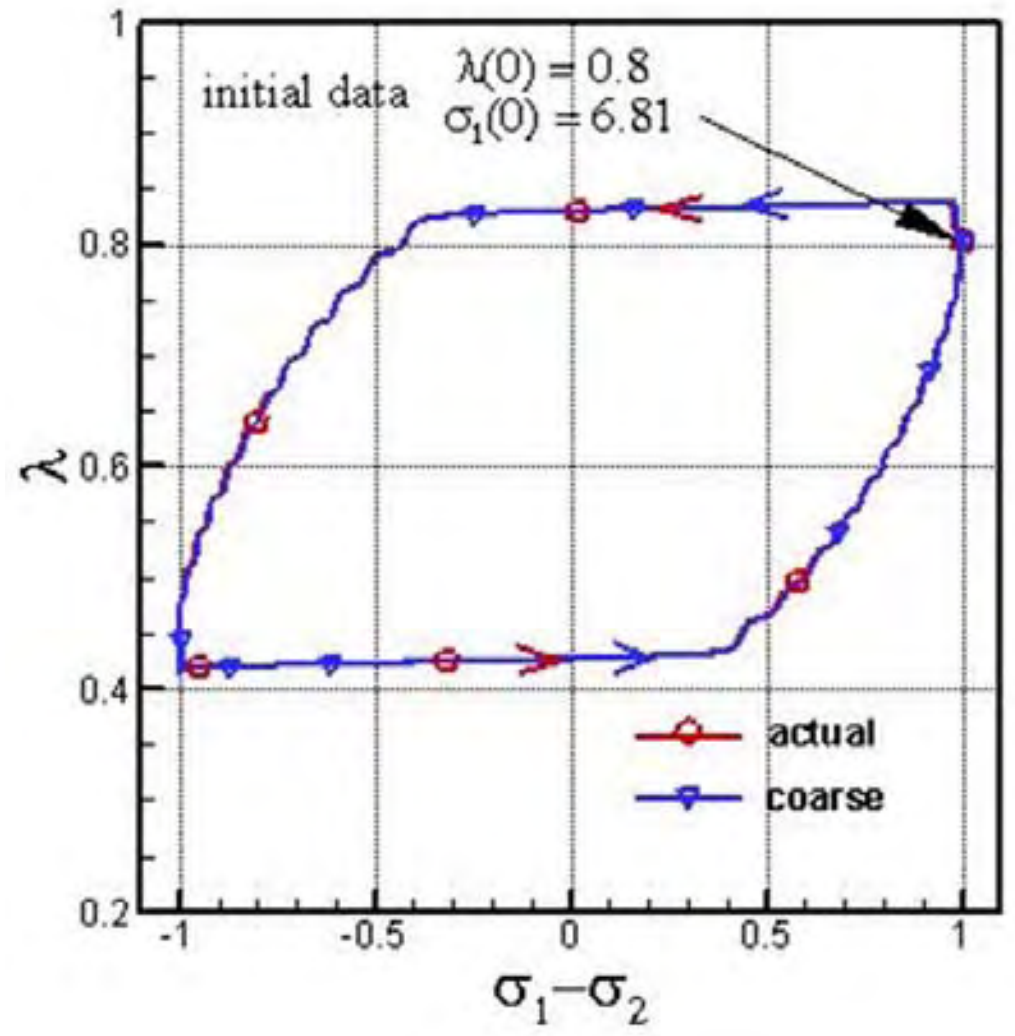

(a)

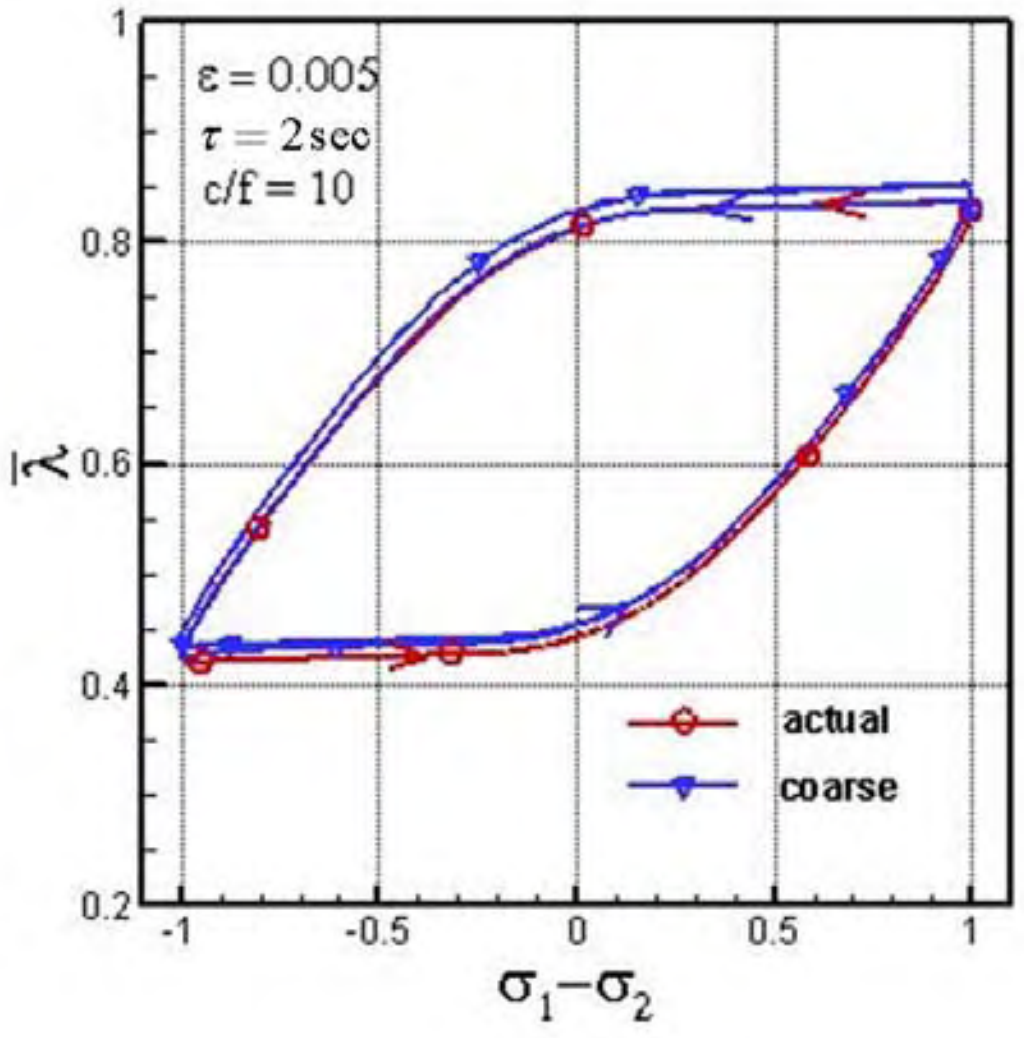

(b)

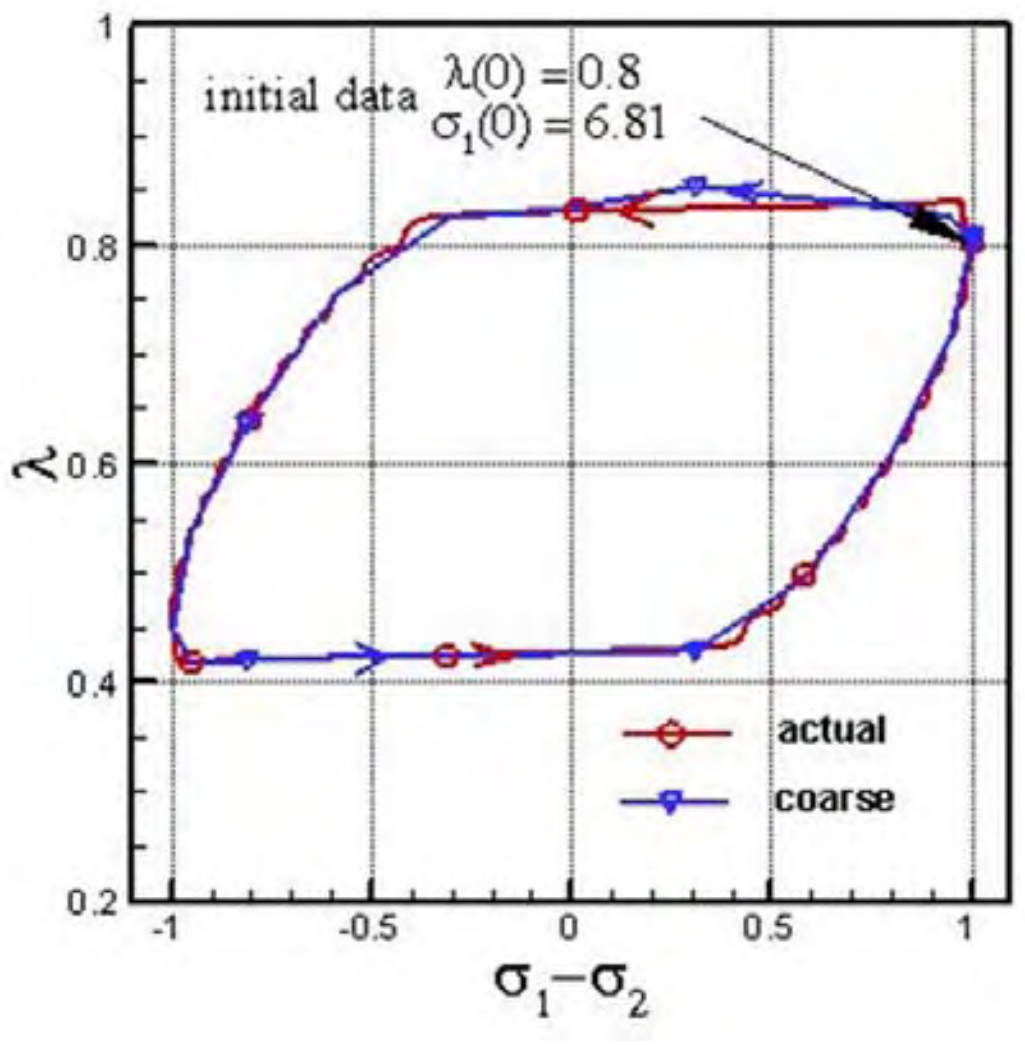

(c)

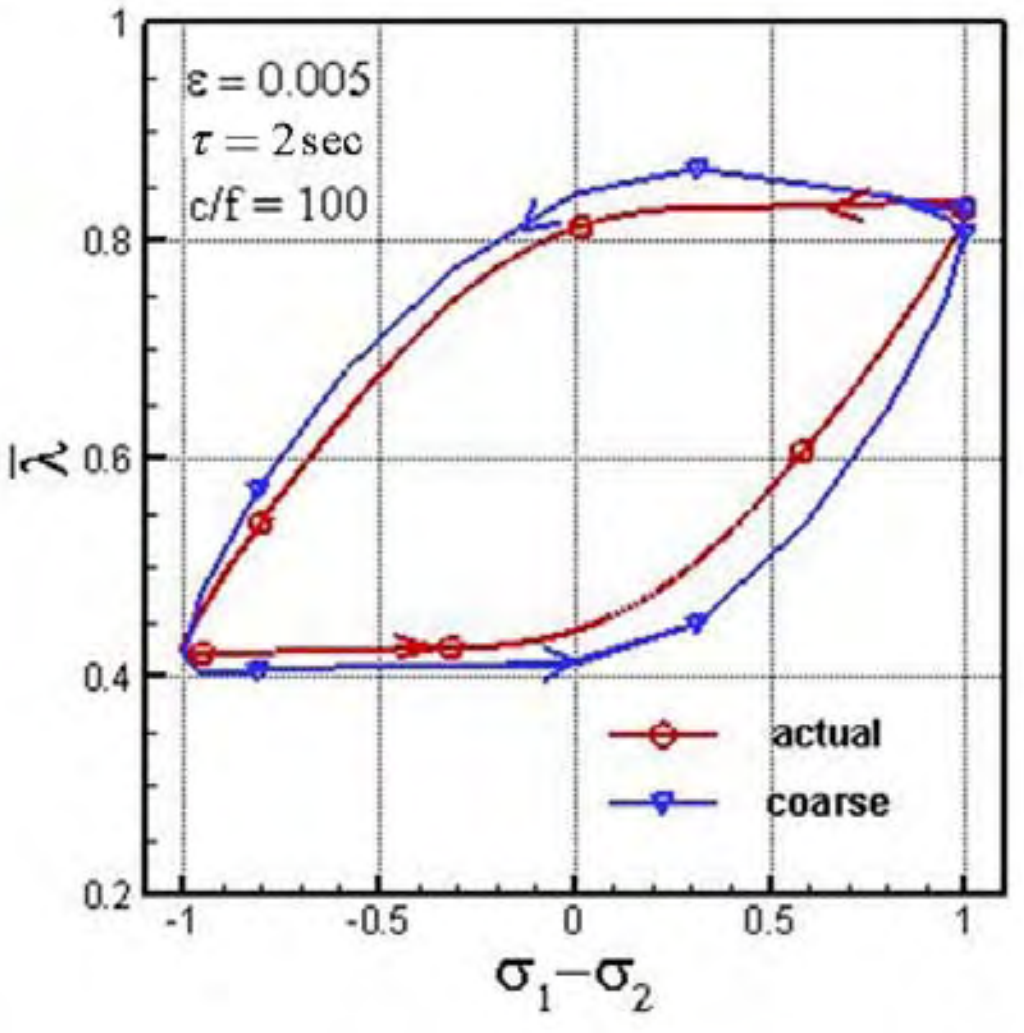

(d)

Figure 4: Comparison of hysteresis for different coarse-to-fine time-step ratios. (a) fine $\lambda\, vs.\, \sigma_1 - \sigma_2$ and (b) averaged $\bar{\lambda}\, vs.\, \sigma_1 - \sigma_2$ for $c/f = 10$: (c) fine $\lambda\, vs.\, \sigma_1 - \sigma_2$ and (d) averaged $\bar{\lambda}\, vs.\, \sigma_1 - \sigma_2$ for $c/f = 100$

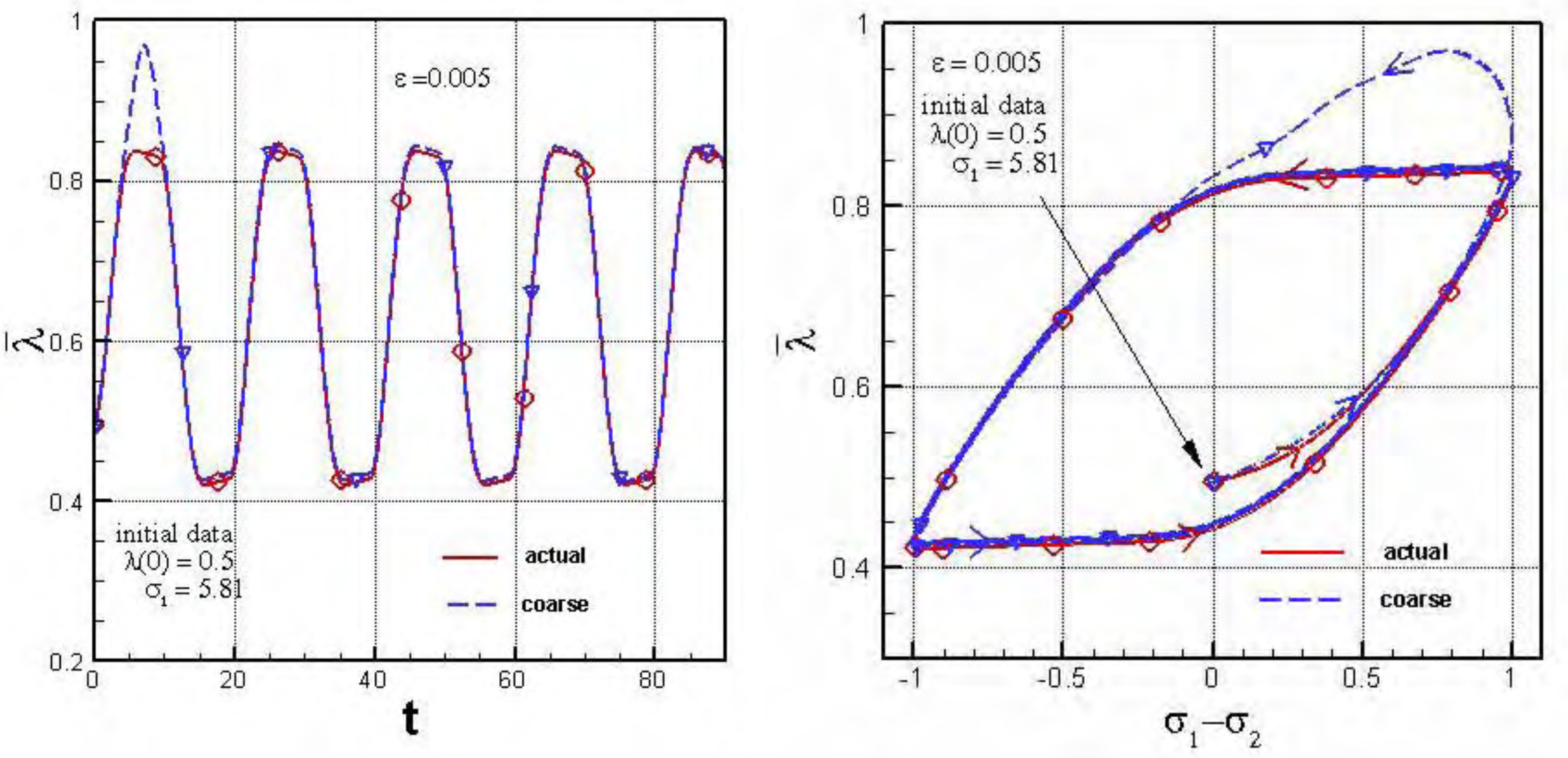

(a) (b)

Figure 5: Discrepancies in coarse response (a) $\bar{\lambda}\, vs.\, t$ , (b) $\bar{\lambda}\, vs.\, \sigma_1 - \sigma_2$

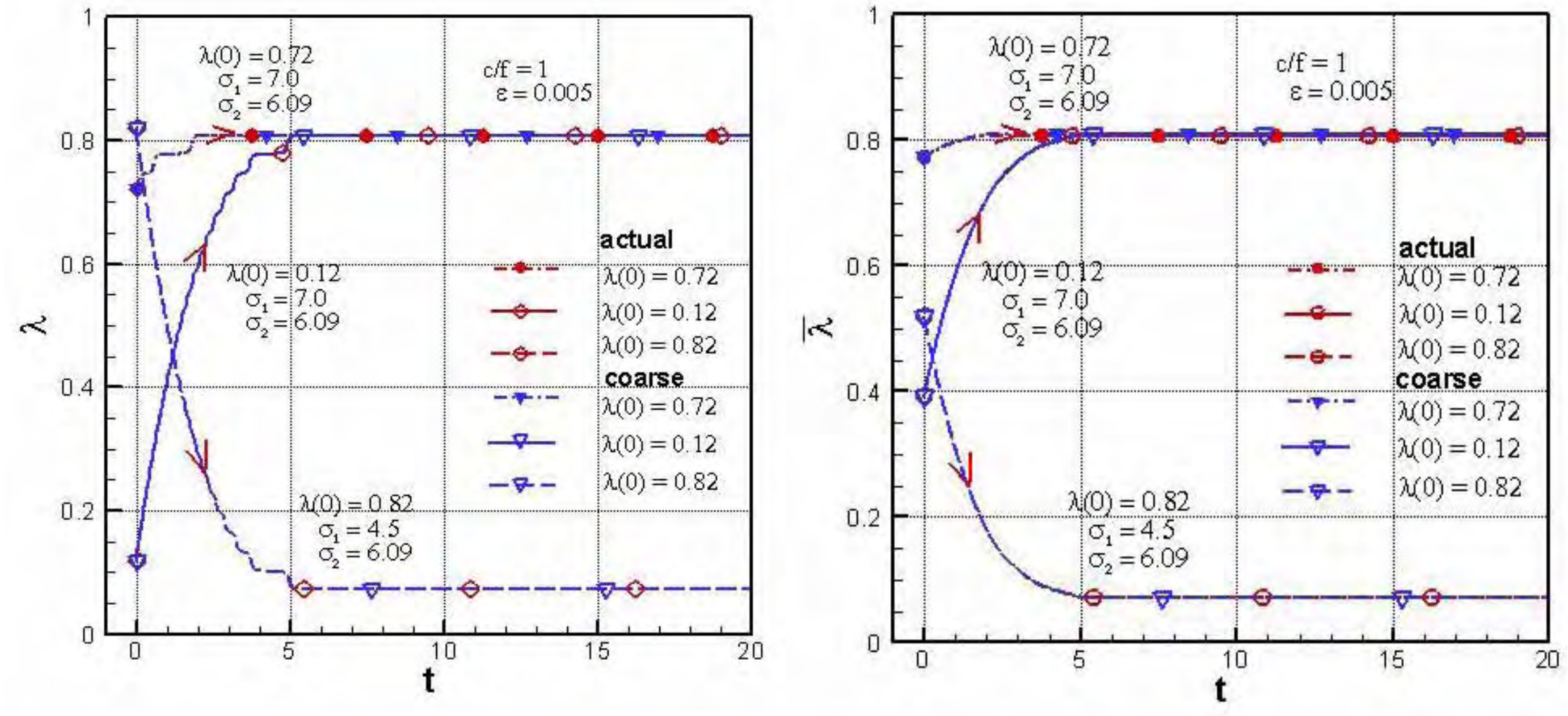

(a) (b)

Figure 6: Creep test: comparison of volume fraction as a function of time
(a) fine $\lambda\, vs.\, t$ , (b) averaged $\bar{\lambda}\, vs.\, t$

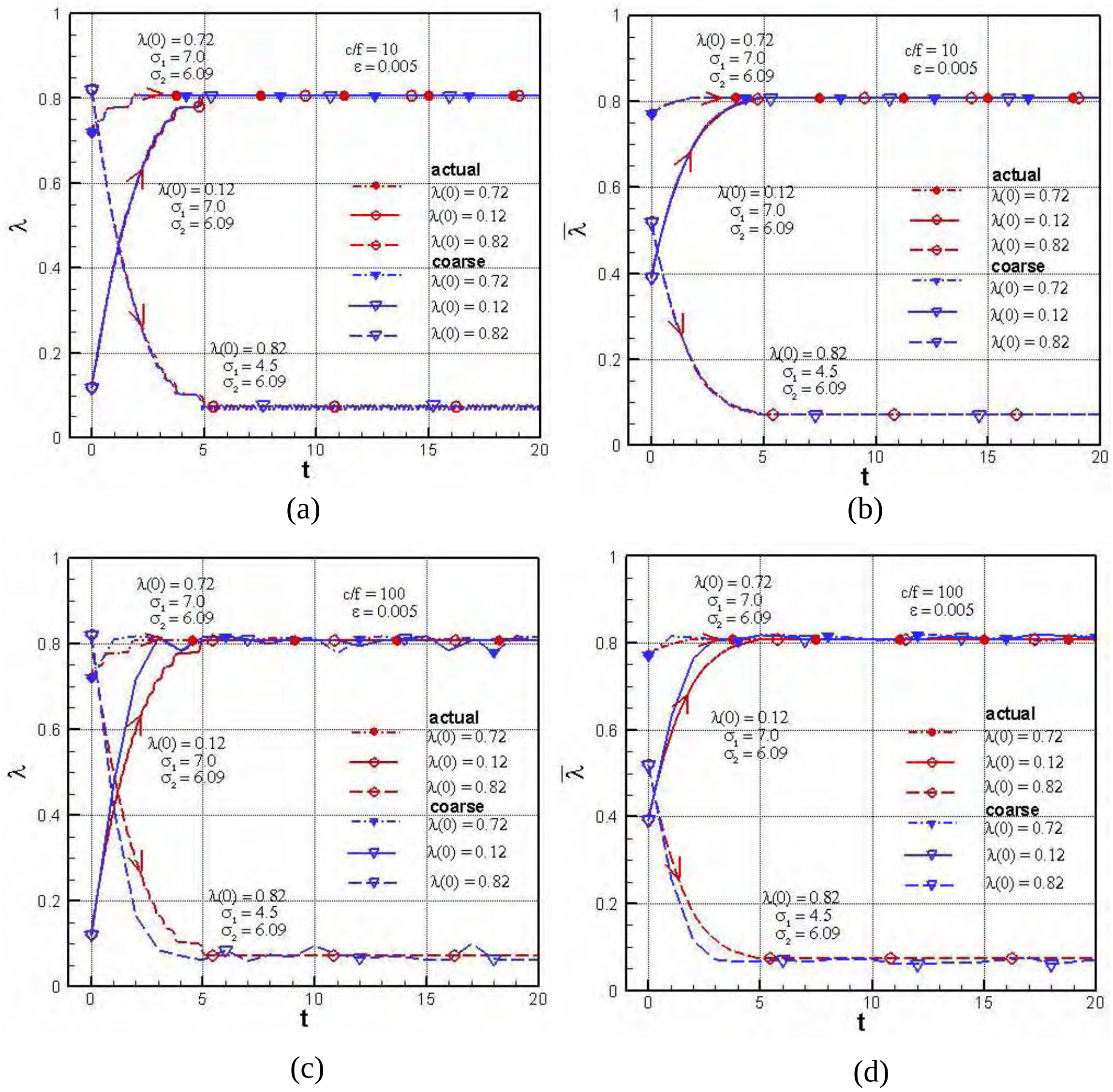

Figure 7: Comparison of volume fraction as a function of time for different coarse-to-fine time-step ratios. (a) fine $\lambda$ *vs.* $\sigma_1 - \sigma_2$ and (b) averaged $\bar{\lambda}$ *vs.* $\sigma_1 - \sigma_2$ for $c/f = 10$. (c) fine $\lambda$ *vs.* $\sigma_1 - \sigma_2$ and (d) averaged $\bar{\lambda}$ *vs.* $\sigma_1 - \sigma_2$ for $c/f = 100$.

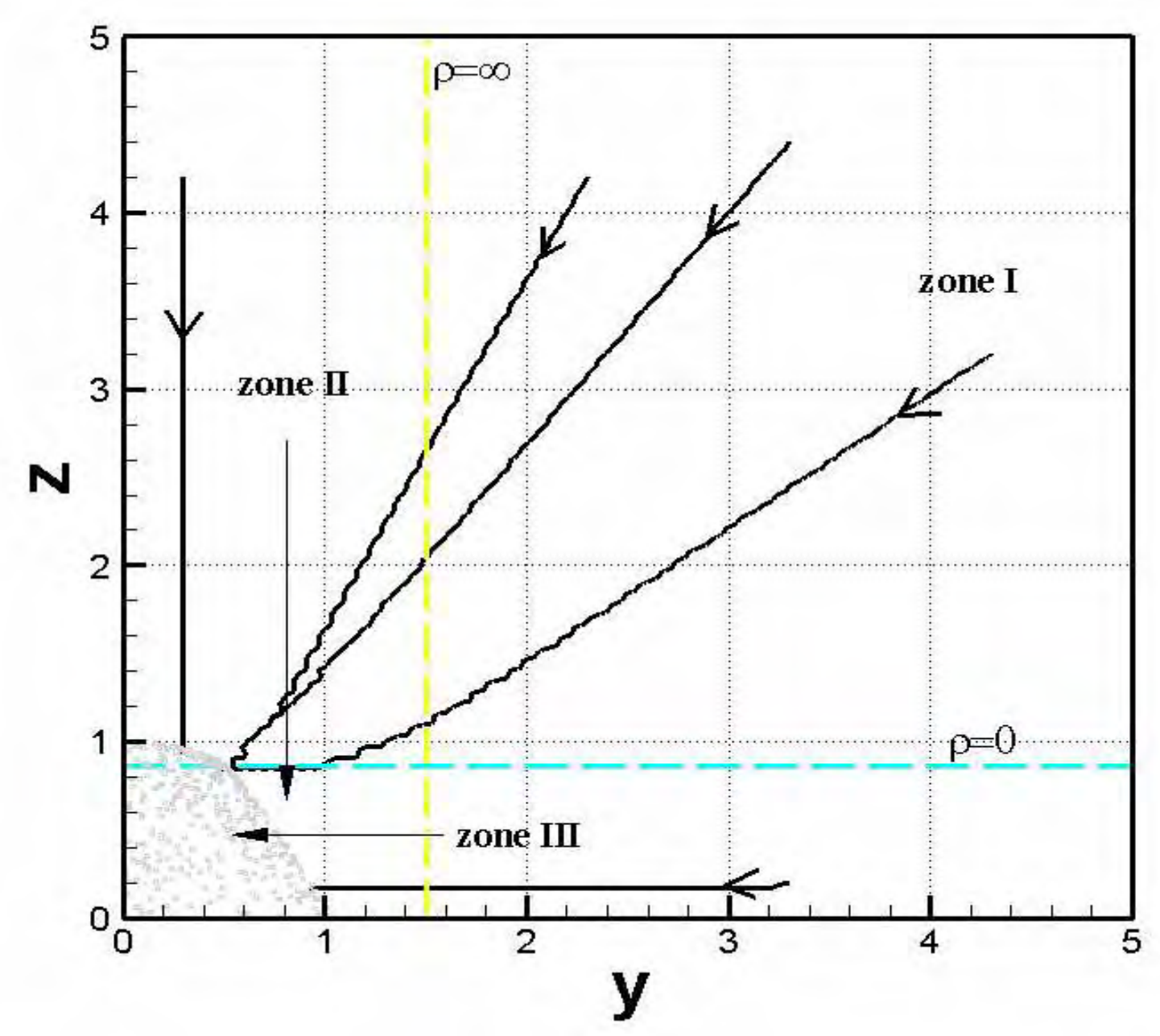

Figure 8: Schematic representation of the fine phase space

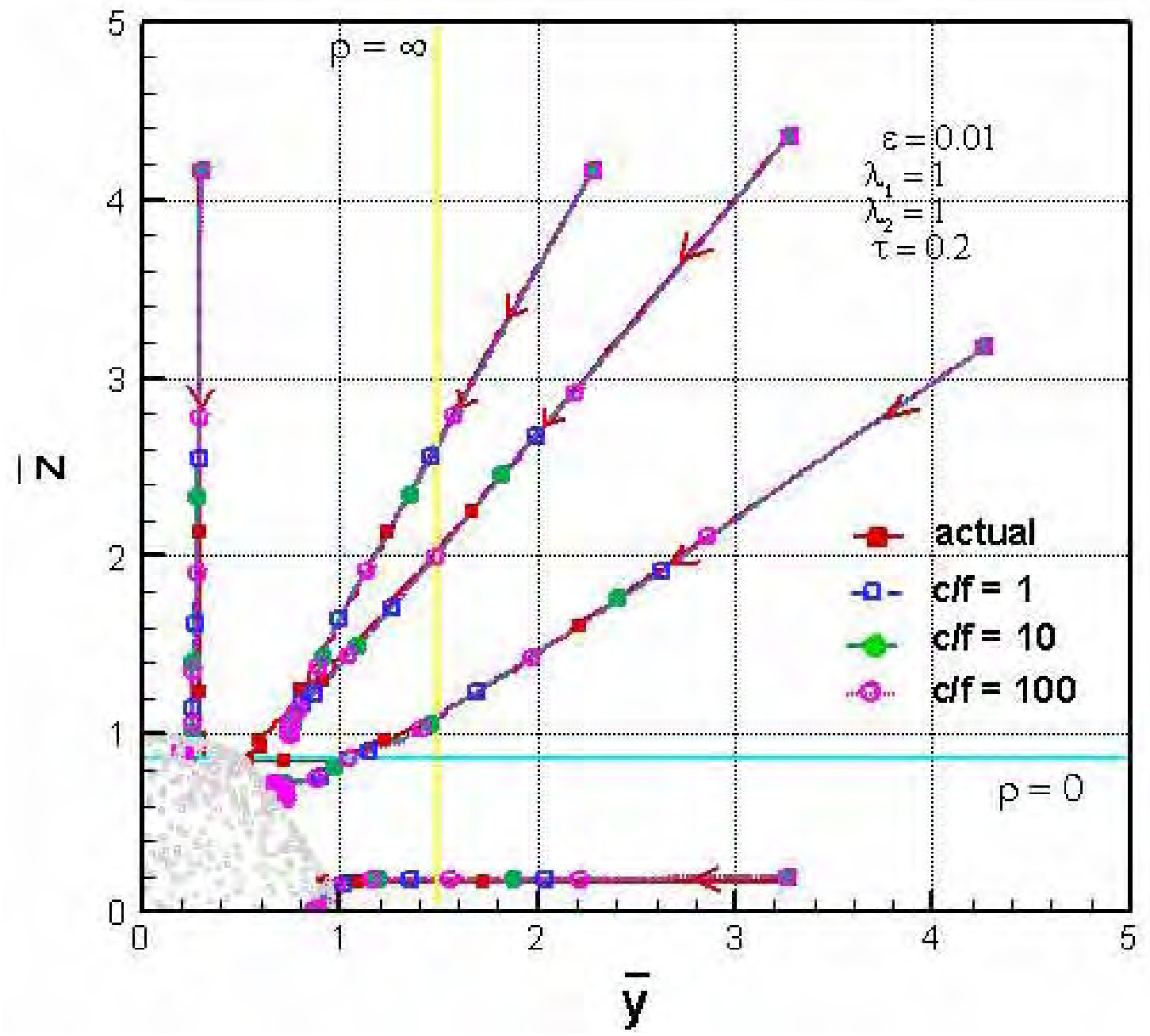

Figure 9: Comparison of coarse and exact response for different coarse-to-fine time-step ratios. $\varepsilon = 0.01$

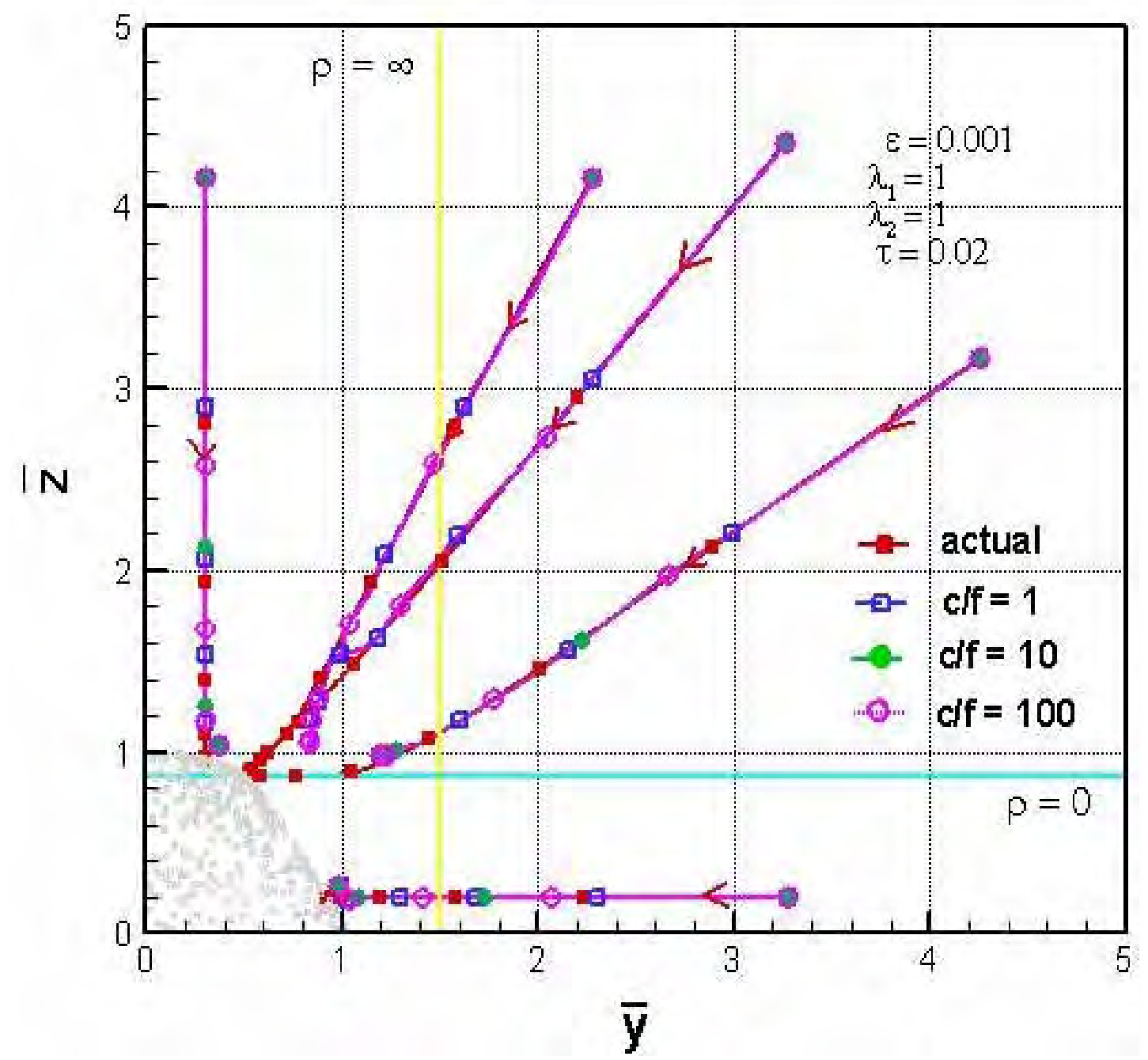

Figure 10: Comparison of coarse and exact response for different coarse-to-fine time-step ratios. $\varepsilon = 0.001$

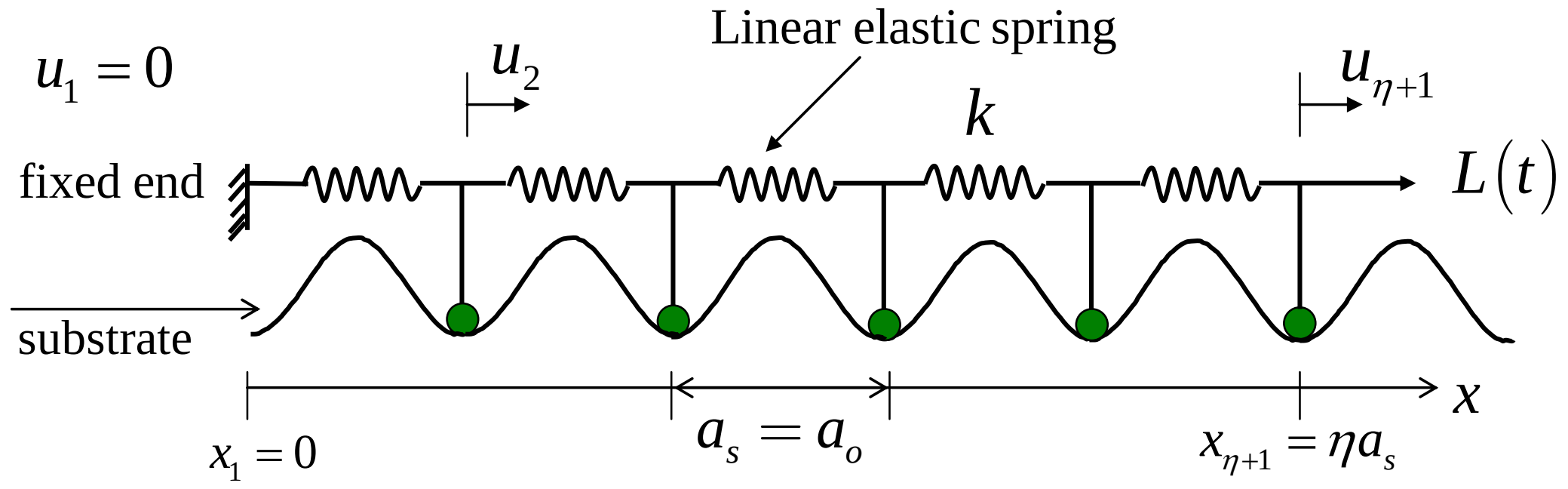

Figure 11: Schematic representation of a chain of atoms using Frenkel-Kontorova model

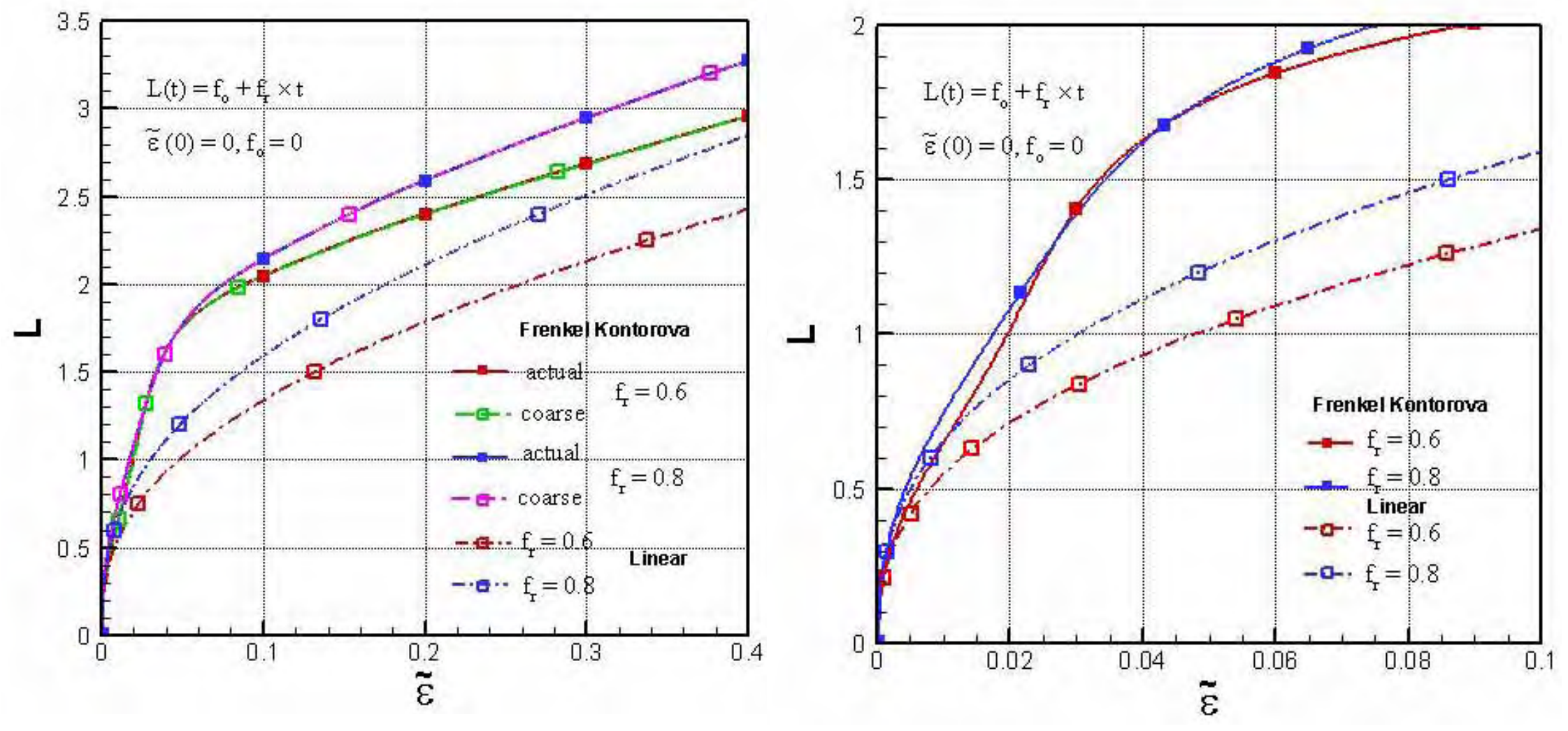

(a) (b)

Figure 12: Comparison of force-strain curve of linear and Frenkel-Kontorova model

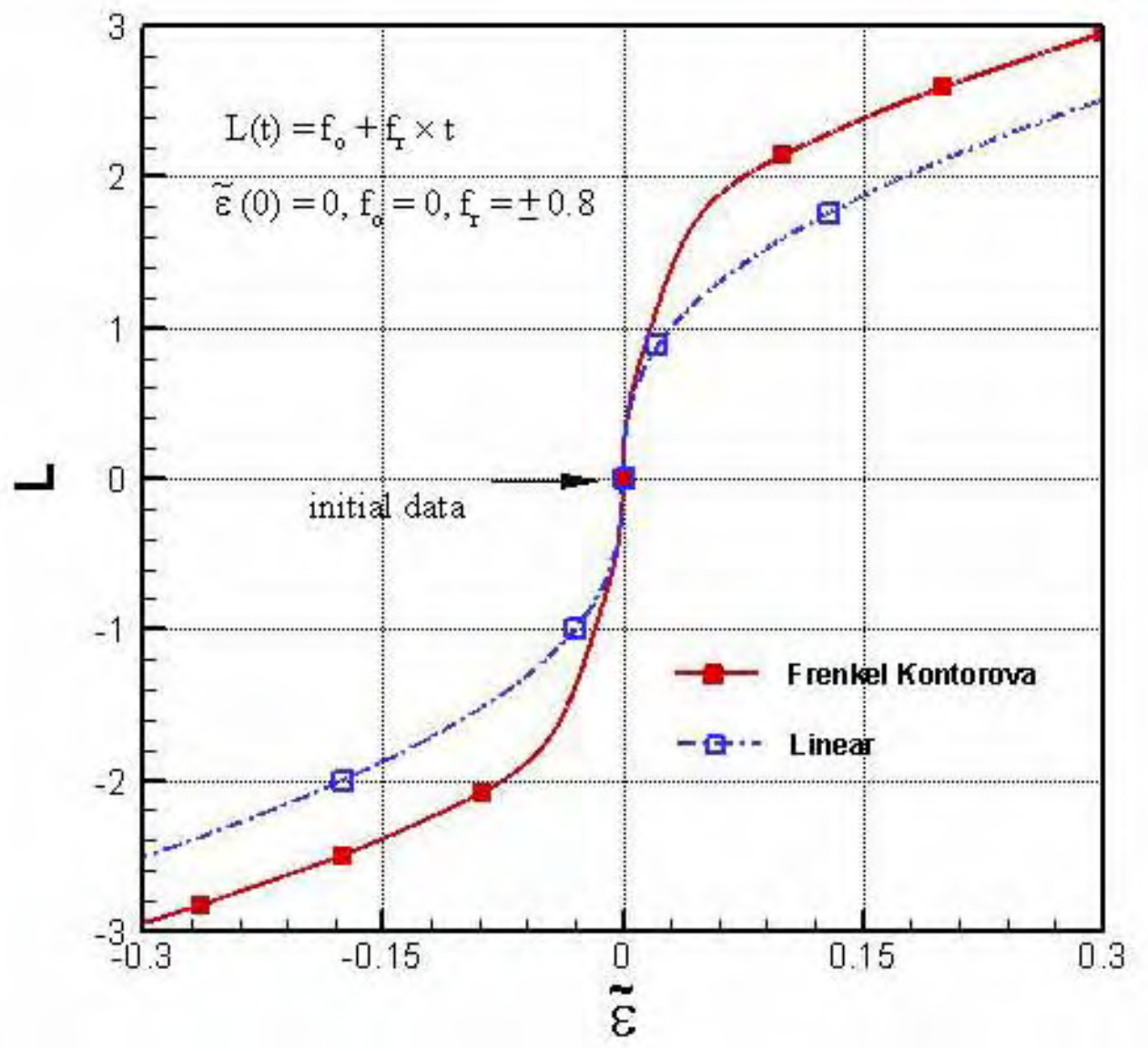

Figure 13: Force-strain curve for loading and unloading from zero ground state

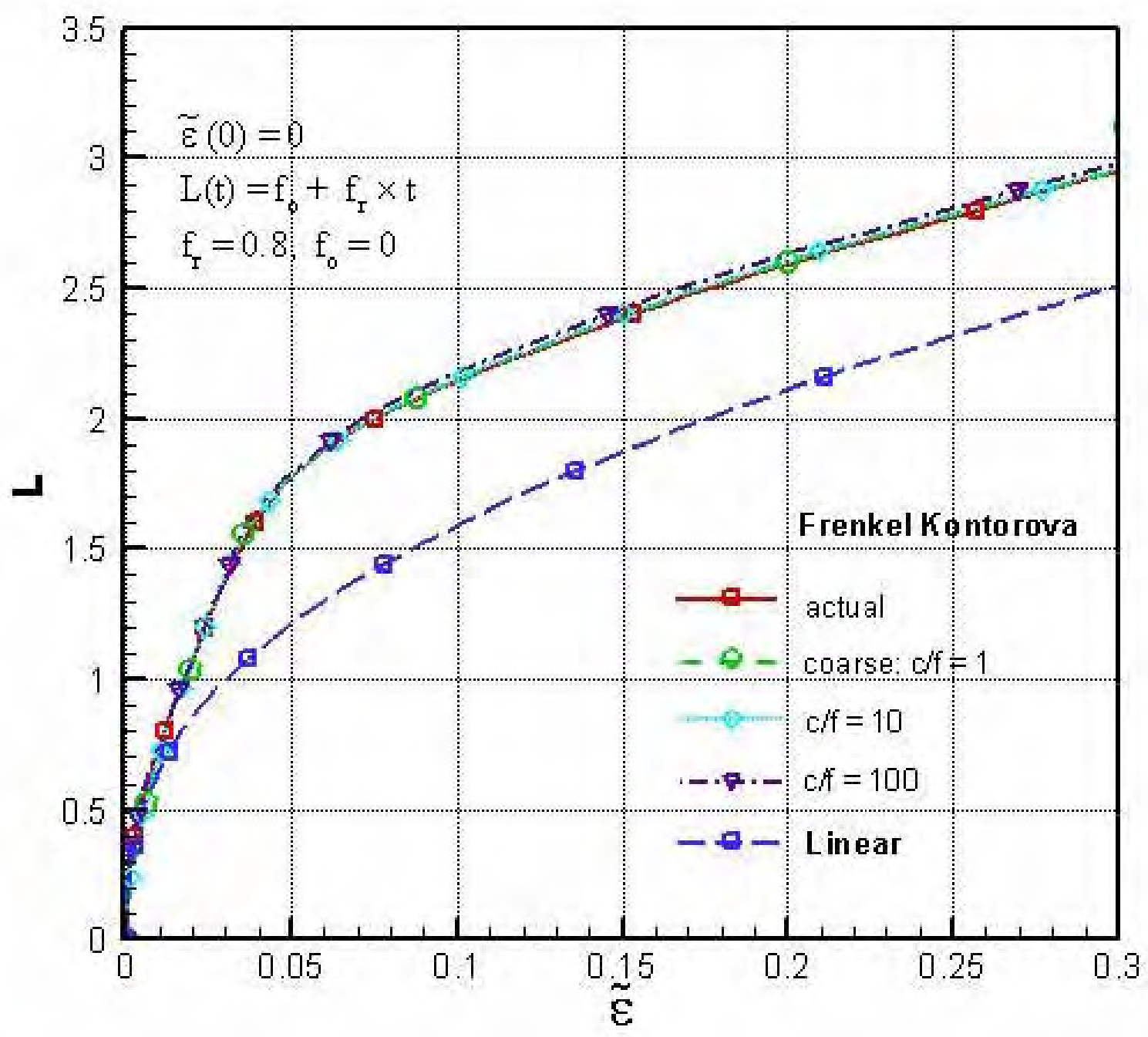

Figure 14: Comparison of coarse and exact response for different coarse-to-fine time-step ratios, for monotonically increasing load example

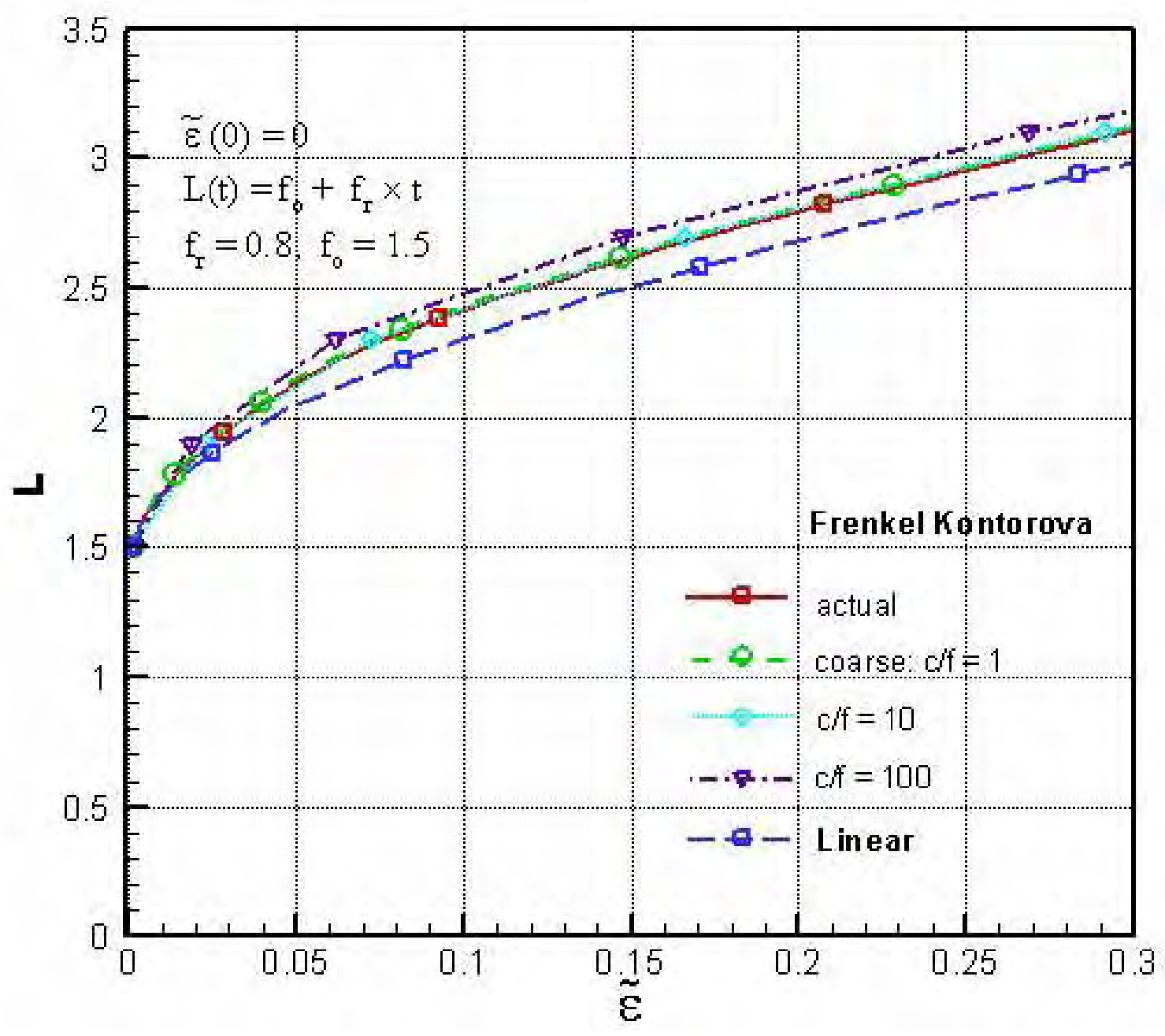

Figure 15: Comparisons of coarse and exact response for different coarse-to-fine time-step ratios, for monotonically increasing load after impulse example

**Appendix:**

**Explicit integration of PDE in the direction of the time-like coarse variable (EI-TC)**

This section gives a general procedure for solving the PDE that arises in the PLIM method by explicit integration in the direction of time-like coarse variable (Sawant, 2005). This procedure was used to compute the locally invariant manifolds for the model problems in Section 3.

Consider an augmented/ non-augmented fine dynamical system (e.g. (6), (15)+(17), (24)+(25), (31)+(35)). Let $\underset{\sim}{\hat{f}} = \left\{ \hat{f}_1, \cdots, \hat{f}_\gamma \right\}$ be the generic state of the fine dynamical system, e.g. $\underset{\sim}{\hat{f}}$ corresponding to (1)-(6) consists of the list $\left( f, f_f \right)$. Let $\underset{\sim}{c} = \left\{ c_1, \cdots, c_\mu \right\}$ be the selected coarse variables. The unknown functions representing the locally invariant manifolds corresponding to $\underset{\sim}{\hat{f}}$ are denoted as $\underset{\sim}{\hat{G}}\left(\underset{\sim}{c}\right) = \underset{\sim}{\hat{f}}$. Let the governing equation for $k^{th}$ unknown function $\hat{G}_k$, similar to (7) described in section 2 (e.g. (18), (26), (36)), be of the form,

$$\sum_{l=1}^{\mu} \frac{\partial \hat{G}_k}{\partial c_l} S_l \left( \underset{\sim}{c}, \underset{\sim}{\hat{G}}\left(\underset{\sim}{c}\right) \right) = H_k \left( \underset{\sim}{\hat{G}}\left(\underset{\sim}{c}\right) \right), \quad f_k^* = \hat{G}_k \left( \left( \underset{\sim}{c}^* \right) \right), \quad k = 1 \, to \, \gamma \tag{a-41}$$

where, the term $S_l$ represents the $rhs$ of the evolution equation for the $l^{th}$ coarse variable, e.g. (11), (19), (27), (37). In the constraint equation (a-41)$_2$, $f_k^*$ is the imposed value for the $k^{th}$ function $\hat{G}_k$ at coarse state $\underset{\sim}{c}^*$.

*Note that,* $\underset{\sim}{\hat{G}}$ *corresponds to all the unknown functions to be solved using* (a-41)*. These should not be confused with* $G, G_f$ *in Section 2. For clarity the carets over all the terms are dropped in rest of this section.*

To solve the governing equations (a-41), we choose a 'time-like' variable from the set $\underset{\sim}{c}$, denoted by $c_t$, say $c_t = c_m,\ 1 \le m \le \mu$. One requirement on this choice is that the *rhs or S* of its evolution equation, is non-zero (at least at all points of the initial hyperplane $c_t = c_t^*$). The remaining set of coarse variables are defined as,

$$\underset{\sim}{x} = \left\{ x_1, \cdots, x_{\mu-1} \right\} \tag{a-42}$$

The discrete approximation for the $k^{th}$ function representing the locally invariant manifolds, $G_k\left(\underset{\sim}{x}; c_t\right)$, is of the following form,

$$G_k\left(\underset{\sim}{x}, c_t\right) \approx \sum_{i}^{n} g_k^i \left(c_t\right) \varphi^i \left(\underset{\sim}{x}\right) \quad \underset{\sim}{x} \in \Omega^{\mu-1}. \tag{a-43}$$

Here $n$ is the total number of discrete nodes in the reduced coarse phase space $\Omega^{\mu-1}$, $x_k$ is the $k^{th}$ independent coarse variable and $g_k^i\left(c_t\right)$ is the value of the function $G_k$ evaluated at node $\underset{\sim}{x}^i$ at 'time' $c_t$. A typical governing equation and constraint equation for $G_k$ can be written as,

$$\frac{\partial G_k}{\partial c_t}\left(\underset{\sim}{x},c_t\right)S_t\left(\underset{\sim}{x},c_t,\underset{\sim}{G}\left(\underset{\sim}{x},c_t\right)\right)+\sum_{l=1}^{\mu-1}\frac{\partial G_k}{\partial x_l}\left(\underset{\sim}{x},c_t\right)S_l\left(\underset{\sim}{x},c_t,\underset{\sim}{G}\left(\underset{\sim}{x},c_t\right)\right)=H_k\left(\underset{\sim}{G}\left(\underset{\sim}{x},c_t\right)\right),\quad k=1\,to\,\gamma \tag{a-44}$$

Where $S_t$ is the non-zero *rhs* of the evolution equation for coarse variable $c_t$. It is important to note that due to the requirement of non-zero $S_t$, one has to choose a different coarse variable as $c_t$ in case we reach a point or curve where $S_t=0$.

The derivative of $G_k$ with respect to $c_t$ is approximated by

$$\frac{\partial G}{\partial c_t}\left(\underset{\sim}{x},c_t\right)\approx\frac{G\left(\underset{\sim}{x},c_t+\Delta c_t\right)-G\left(\underset{\sim}{x},c_t\right)}{\Delta c_t}, \tag{a-45}$$

and this discretization is substituted in the governing equation (a-44) to compute the solution. For this approach, consistent data is imposed at an initial step $c_t^*$ over the coarse domain $\Omega^{\mu-1}$ as,

$$f_k^*=G_k\left(\underset{\sim}{x},c_t^*\right),\quad k=1\,to\,\gamma. \tag{a-46}$$

The following figure illustrates the above discretization for the model problem 2 of 2-D gradient system with wiggly energy. The augmented fine system consists of four fine variables $y,y_f,z,z_f$, i.e. $\gamma=4$. The coarse variables are $\underset{\sim}{c}=\left\{\bar{y},\bar{z}\right\}$, i.e. $\mu=2$.

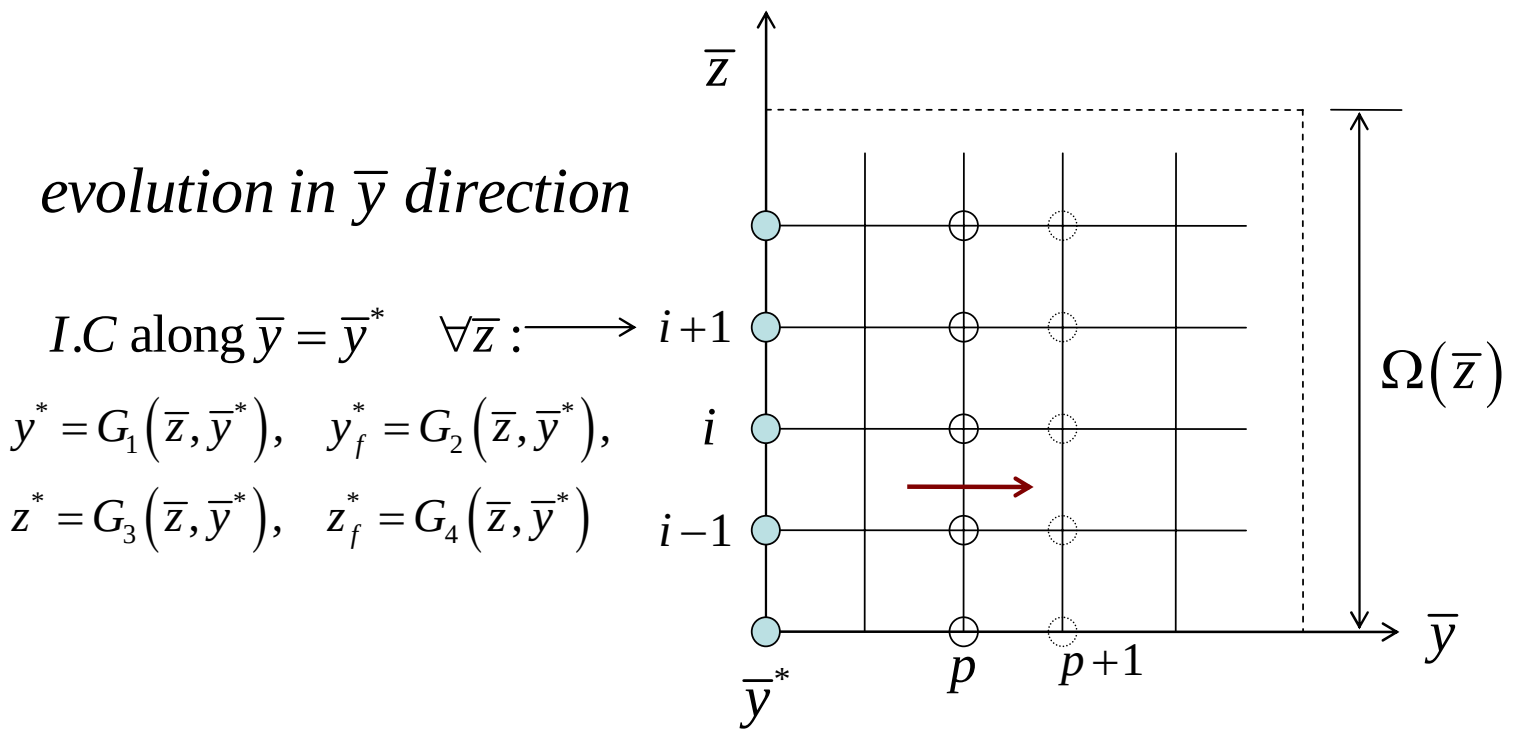

Schematic diagram for evolution in the direction of time-like coarse variable (EI-TC)

Here we choose say $\bar{y}$ as the time-like variable and the unknown functions $G_k\left(\bar{y},\bar{z}\right)$ are discretized in $\bar{z}$ direction using finite elements. The solutions to $G_k\left(\bar{y},\bar{z}\right)$ are computed by marching in the $\bar{y}$ direction, where $\bar{y} = \bar{y}^*$ is the initial curve over which consistent data is imposed that serve as initial conditions for the 'time' marching scheme.

From (a-44), (a-45), (a-46), we construct the following quadratic functional that is defined in terms of the $L_2$ norms of the equation residuals:

$$I : \Xi^\gamma \to \mathbb{R}$$

$$I\left(\underset{\sim}{G}^{p+1}\right) = \int_{\Omega^{\mu-1}} \left\{ \sum_{k=1}^{\gamma} \left( \begin{array}{l} \dfrac{G_k^{p+1}\left(\underset{\sim}{x}\right) - G_k^p\left(\underset{\sim}{x}\right)}{\Delta c_t} S_t\left(\underset{\sim}{x}, c_t, \underset{\sim}{G}^p\left(\underset{\sim}{x}\right)\right) \\ + \displaystyle\sum_{l=1}^{\mu-1} \dfrac{\partial G_k^{p+1}\left(\underset{\sim}{x}\right)}{\partial x_l} S_l\left(\underset{\sim}{x}, c_t, \underset{\sim}{G}^p\left(\underset{\sim}{x}\right)\right) - H_k\left(\underset{\sim}{G}^p\left(\underset{\sim}{x}\right)\right) \end{array} \right)^2 \right\} d\underset{\sim}{x}, \tag{a-47}$$

where $G_k^p\left(\underset{\sim}{x}\right), G_k^{p+1}\left(\underset{\sim}{x}\right)$ are the values of the function $G_k\left(\underset{\sim}{x}\right)$ evaluated at $p^{th}$ and $p+1^{th}$ step in the time-like direction and the terms $S, H$ are computed from the $\underset{\sim}{G}^p\left(\underset{\sim}{x}\right)$ at $p^{th}$ step. Thus the EI-TC method is a linearized numerical scheme with consistent mixture of forward & backward Euler schemes.

A necessary condition that $\underset{\sim}{G}^{p+1} \in \Xi^\gamma$ is a minimizer of the functional $I$ is that its first variation vanishes at $\underset{\sim}{G}^{p+1}$ for all admissible $\underset{\sim}{Q} \in \Xi^\gamma$. This leads to the least-squares variational formulation: Find $\underset{\sim}{G}^{p+1} = \left\{G_1^{p+1}, \cdots, G_\gamma^{p+1}\right\} \in \Xi^\gamma$ such that

$$\lim_{\delta \to 0} \frac{d}{d\delta} I\left(G_k^{p+1} + \delta Q_k\right) =$$

$$\int_{\Omega^{\mu-1}} \left[ \begin{array}{l} \left( \dfrac{G_k^{p+1}\left(\underset{\sim}{x}\right) - G_k^p\left(\underset{\sim}{x}\right)}{\Delta c_t} S_t\left(\underset{\sim}{x}, c_t\right) + \displaystyle\sum_{l=1}^{\mu-1} \dfrac{\partial G_k^{p+1}\left(\underset{\sim}{x}\right)}{\partial x_l} S_l\left(\underset{\sim}{x}, c_t\right) - H_k\left(\underset{\sim}{G}^p\left(\underset{\sim}{x}\right)\right) \right) \cdot \\ \left( \dfrac{Q_k\left(\underset{\sim}{x}\right)}{\Delta c_t} S_t\left(\underset{\sim}{x}, c_t\right) + \displaystyle\sum_{l=1}^{\mu-1} \dfrac{\partial Q_k\left(\underset{\sim}{x}\right)}{\partial x_l} S_l\left(\underset{\sim}{x}, c_t\right) \right) \end{array} \right] d\underset{\sim}{x} = 0 \tag{a-48}$$

$$\forall Q \in \Xi^\gamma$$

Thus the simplified coupled system to be solved is of the following form,

$$\left[M_{kk}\right]^{p}\left[g_{k}\right]^{p+1}=\left[\alpha_{k}\right]^{p}\quad k=1\,to\,\gamma$$

$$\left[M_{kk}\right]_{ij}^{p}=\int_{\Omega^{\mu-1}}\left[\begin{array}{l}\left(\varphi^{i}\,\dfrac{S_{t}\left(\underset{\sim}{x},c_{t}\right)}{\Delta c_{t}}+\sum_{l=1}^{\mu-1}\dfrac{\partial\varphi^{i}}{\partial x_{l}}S_{l}\left(\underset{\sim}{x},c_{t}\right)\right).\\ \left(\varphi^{j}\,\dfrac{S_{t}\left(\underset{\sim}{x},c_{t}\right)}{\Delta c_{t}}+\sum_{l=1}^{\mu-1}\dfrac{\partial\varphi^{j}}{\partial x_{l}}S_{l}\left(\underset{\sim}{x},c_{t}\right)\right)\end{array}\right]d\underset{\sim}{x}$$

$$\left[\alpha_{k}\right]_{i}^{p}=\int_{\Omega^{\mu-1}}\left[\left(\varphi^{i}\,\frac{S_{t}\left(\underset{\sim}{x},c_{t}\right)}{\Delta c_{t}}+\sum_{l=1}^{\mu-1}\frac{\partial\varphi^{i}}{\partial x_{l}}S_{l}\left(\underset{\sim}{x},c_{t}\right)\right).\left(G_{k}^{p}\,\frac{S_{t}\left(\underset{\sim}{x},c_{t}\right)}{\Delta c_{t}}+H_{k}\left(\underset{\sim}{G}^{p}\right)\right)\right]d\underset{\sim}{x}$$

$$\left[g_{k}\right]^{p+1}=\left\{g_{k}^{1},\cdots,g_{k}^{n}\right\}^{p+1}\quad i,\,j=1\,to\,n.$$

(a-49)